\documentclass[11pt]{amsart}
\usepackage{amsxtra}
\usepackage{amssymb}
\theoremstyle{plain}

\newcommand{\cleqn}{\setcounter{equation}{0}}
\newcommand{\clth}{\setcounter{theorem}{0}}
\newcommand {\sectionnew}[1]{\section{#1}\cleqn\clth}
\newcommand{\nn}{\hfill\nonumber}
\newtheorem{theorem}{Theorem}[section]
\newtheorem{lemma}[theorem]{Lemma}
\newtheorem{definition-theorem}[theorem]{Definition-Theorem}
\newtheorem{proposition}[theorem]{Proposition}
\newtheorem{corollary}[theorem]{Corollary}
\newtheorem{definition}[theorem]{Definition}
\newtheorem{example}[theorem]{Example}
\newtheorem{remark}[theorem]{Remark}
\newtheorem{conjecture}[theorem]{Conjecture}
\newtheorem{notation}[theorem]{Notation}
\newcommand \bth[1] { \begin{theorem}\label{t#1} }
\newcommand \ble[1] { \begin{lemma}\label{l#1} }

\newcommand \bpr[1] { \begin{proposition}\label{p#1} }
\newcommand \bco[1] { \begin{corollary}\label{c#1} }
\newcommand \bde[1] { \begin{definition}\label{d#1}\rm }
\newcommand \bex[1] { \begin{example}\label{e#1}\rm }
\newcommand \bre[1] { \begin{remark}\label{r#1}\rm }
\newcommand \bcj[1] { \begin{conjecture}\label{j#1}\rm }

\newcommand \bnota[1] { \begin{notation}\label{n#1}\rm }
\renewcommand {\eth} { \end{theorem} }
\newcommand {\ele} { \end{lemma} }

\newcommand {\epr} { \end{proposition} }
\newcommand {\eco} { \end{corollary} }
\newcommand {\ede} { \end{definition} }
\newcommand {\eex} { \end{example} }
\newcommand {\ere} { \end{remark} }
\newcommand {\ecj} { \end{conjecture} }

\newcommand {\enota} { \end{notation} }
\newcommand \thref[1]{Theorem \ref{t#1}}
\newcommand \leref[1]{Lemma \ref{l#1}}
\newcommand \prref[1]{Proposition \ref{p#1}}

\newcommand \deref[1]{Definition \ref{d#1}}

\newcommand \lb[1]{\label{#1}}

\def \Rset {{\mathbb R}}         

\def \KK {{\mathbb K}}
\def \AA {{\mathbb A}}
\def \Zset {{\mathbb Z}}
\def \Nset {{\mathbb N}}
\def \Qset {{\mathbb Q}}

\def \Tset {{\mathbb T}}

\def \B  {{\mathcal{B}}}               

\def \CD {{\mathcal{CD}}}
\def \LP {{\mathcal{LP}}}
\def \QQ {{\mathcal{Q}}}
\def \PP {{\mathcal{P}}}

\def \UU {{\mathcal{U}}}
\def \RR {{\mathcal{R}}}
\def \SS {{\mathcal{S}}}

\def \TT {{\mathcal{T}}} 

\def \ib {{\bf{i}}}
\def \De {\Delta}   
\def \de {\delta}
\def \al {\alpha}
\def \be {\beta}
\def \vpi {\varpi}
\def \la {\lambda}

\def \kap {\kappa}
\def \om {\omega}
\def \Om {\Omega}
\def \ga {\gamma}
\def \de {\delta}

\def \sig {\sigma}

\def \vp {\varphi}
\def \ep {\epsilon}
\def \sig{\sigma}
\def \qb { {\bf{q}} }

\def \mt  {\mapsto}

\def \hra {\hookrightarrow}

\def \rcor {\rangle}
\def \lcor {\langle}
\def \o  {\otimes}
\def \ol {\overline}

\def \wh {\widehat}

\def \id { {\mathrm{id}} }

\def \g  {\mathfrak{g}}   

\def \sl {\mathfrak{sl}}

\def \n  {\mathfrak{n}}

\def \b  {\mathfrak{b}}

\def \sl {\mathfrak{sl}}

\DeclareMathOperator \congto  {{ \stackrel{\cong}{\to} }}
\DeclareMathOperator \Span { {\mathrm{Span}} }
\DeclareMathOperator \Aut { {\mathrm{Aut}} }

\DeclareMathOperator \GKdim {{\mathrm{GK \, dim}}}

\DeclareMathOperator \lt  { {\mathrm{lt}} }

\DeclareMathOperator \Fract { {\mathrm{Fract}} }
\DeclareMathOperator \TSpec { {\mathbb T}^r \mbox{-} {\mathrm{Spec}}}

\renewcommand \max { {\mathrm{max}} }
\newcommand \Spec { {\mathrm{Spec}} }

\begin{document}
\title[Quantum Schubert cells]
{Quantum Schubert cells via \\ representation theory and ring theory}
\author[Joel Geiger]{Joel Geiger}
\address{
Department of Mathematics \\
Louisiana State University \\
Baton Rouge, LA 70803
U.S.A.
}
\email{jgeige1@math.lsu.edu}
\author[Milen Yakimov]{Milen Yakimov}
\email{yakimov@math.lsu.edu}
\date{}
\begin{abstract}
We resolve two questions of Cauchon and M\'eriaux on the spectra 
of the quantum Schubert cell algebras $\UU^-[w]$. 
The treatment of the first one unifies two very different approaches 
to $\Spec \UU^-[w]$, a ring theoretic 
one via deleting derivations and a representation theoretic one via 
Demazure modules. The outcome is that now one can combine 
the strengths of both methods. As an application we 
solve the containment problem for the Cauchon--M\'eriaux 
classification of torus invariant prime ideals of $\UU^-[w]$.
Furthermore,
we construct explicit models in terms of quantum 
minors for the Cauchon quantum affine 
space algebras constructed via the procedure of 
deleting derivations from all quantum Schubert cell algebras 
$\UU^-[w]$. Finally, our methods also give a new, independent proof of 
the Cauchon--M\'eriaux classification.
\end{abstract}
\maketitle
\sectionnew{Introduction}
\lb{intro}
The study of the spectra of quantum groups for generic deformation 
parameters was initiated twenty years ago by 
Joseph \cite{J0,J} and Hodges--Levasseur--Toro \cite{HLT} who obtained 
a number of important results on them. One of the long-term goals 
was to understand these spectra geometrically in terms of symplectic foliations 
in an attempt to extend the orbit method \cite{Dix} 
to more general classes of algebras and Poisson manifolds.
This grew into a very active area of studying 
the ring theoretic properties of quantum analogs of universal enveloping algebras 
of solvable Lie algebras. The quantum Schubert cell algebras,
defined by De Concini--Kac--Procesi \cite{DKP} and Lusztig \cite{L},
comprise one of the major families of algebras in this area.
There is one such algebra $\UU^-[w]$ 
for every simple Lie algebra $\g$ and an element $w$ of the Weyl group $W$ of $\g$.
It is a subalgebra of the quantized universal enveloping algebra $\UU_q(\g)$ 
and a deformation of the universal enveloping 
algebra $\UU(\n_- \cap w(\n_+))$, where $\n_\pm$ are the nilradicals
of a pair of opposite Borel subalgebras $\b_\pm$ of $\g$. 
From another perspective, the algebra $\UU^-[w]$ is a deformation
of the coordinate ring of the Schubert cell corresponding to $w$ 
of the full flag variety of $\g$, equipped with the standard Poisson 
structure \cite{GY}. These 
algebras played important roles in many different contexts
in recent years such as the 
study of coideal subalgebras of $\UU_q(\b_-)$ and $\UU_q(\g)$  
\cite{HS,HK} and quantum cluster algebras \cite{GLS}.

There are two very different approaches to the study 
of the spectra of $\UU^-[w]$. One is purely ring theoretic and 
is based on the Cauchon procedure of deleting derivations 
\cite{Ca}. The second is a representation theoretic one
and builds on the above mentioned methods of Joseph, Hodges, 
Levasseur, and Toro \cite{J,HLT}. Each of these methods has 
a number of advantages over the other, and relating them 
was an important open problem with many potential applications.
Previously there were no 
connections between them even for special cases of 
the algebras $\UU^-[w]$, such as the algebras of quantum matrices.

In this paper we unify the ring theoretic
and the representation theoretic approaches 
to the study of $\Spec \UU^-[w]$. Furthermore, we resolve several 
other open problems on the deleting derivation procedure and the 
spectra of $\UU^-[w]$, two being questions 
posed by Cauchon and M\'eriaux \cite{MC}.
Before we proceed with the statements of these
results, we need to introduce some additional background.

There is a canonical action of the torus $\Tset^r = (\KK^*)^{\times r}$
on $\UU^-[w]$ by algebra automorphisms, where $\KK$ is the base field 
and $r$ is the rank of $\g$. By a general stratification result of Goodearl 
and Letzter \cite{GL}, one has a partition 
\[
\Spec \UU^-[w] = \bigsqcup_{I \in \TSpec \UU^-[w]} 
\Spec_I \UU^-[w].
\]
Here $\TSpec \UU^-[w]$ denotes the set of $\Tset^r$-invariant 
prime ideals. By two general results of \cite{GL} $\TSpec \UU^-[w]$
is finite and each stratum 
\[
\Spec_I \UU^-[w] = \{ L \in \Spec \UU^-[w] \mid 
\cap_{t \in \Tset^r} t \cdot L = I \}
\]
is homeomorphic to the spectrum of a (commutative) Laurent 
polynomial ring. The problem of the description of the Zariski 
topology of $\Spec \UU^-[w]$, however, is wide open. 

The Cauchon method of deleting derivations is a multi-stage
recursive procedure \cite{Ca} beginning with an iterated 
Ore extension $A$ of length $l$ (of a certain general type) equipped with a compatible 
$\Tset^r$-action and ending with a quantum affine space 
algebra $\ol{A}$ with a $\Tset^r$-action. Cauchon constructed in \cite{Ca} a set-theoretic 
embedding of $\Spec A$ into $\Spec \ol{A}$. It restricts 
to a set-theoretic embedding 
$\TSpec A \hra \TSpec \ol{A}$. The 
$\Tset^r$-invariant prime ideals of $A$ are then 
parametrized by some of the subsets of $[1,l]$, 
called Cauchon diagrams. 
The $\Tset^r$-prime ideal of $A$ corresponding to a 
Cauchon diagram $D \subseteq [1,l]$ will be denoted by $J_D$.
The problem of determining 
which subsets of $[1,l]$ arise in this way (i.e., are Cauchon 
diagrams), is the essence of the method and is very 
difficult for each particular class of algebras. 
It was solved for the algebras of quantum matrices 
by Cauchon \cite{Ca} and for
all algebras $\UU^-[w]$ by Cauchon and M\'eriaux \cite{MC}.
To state the latter result, we denote the set of simple roots of $\g$ by $\Pi$ 
and the corresponding simple reflections of $W$ by 
$s_\al$, $\al \in \Pi$. A word 
$\ib=(\al_1, \ldots, \al_l)$ in the alphabet $\Pi$ 
will be called a reduced word for $w$ if 
$s_{\al_1} \ldots s_{\al_l}$ is a 
reduced expression of $w$. Each reduced word $\ib$ for $w$ 
gives rise to 
a presentation of $\UU^-[w]$ as an iterated Ore extension of length $l$.
The subsets of 
$[1,l]$ are index sets for the subwords of $\ib$
by the assignment
$\{ j_1 < \ldots < j_n \} \mt (\al_{j_1}, \ldots, \al_{j_n})$.
We will denote by $\leq$ the (strong) Bruhat order on $W$ 
and set $W^{\leq w} = \{ y \in W \mid y \leq w \}$.
For each $y \in W^{\leq w}$
there exists a unique left positive subword of $\ib$ corresponding to $y$
(see \S 2.2 for its definition and details on Weyl group combinatorics). 
Its index set will be denoted 
by $\LP_\ib(y)$. The Cauchon--M\'eriaux classification theorem 
states the following:

{\em{For all Weyl group elements $w \in W$ and reduced 
words $\ib$ for $w$, consider the presentation
of $\UU^-[w]$ as an iterated Ore extension corresponding to $\ib$.
The Cauchon diagrams of the 
$\Tset^r$-prime ideals of $\UU^-[w]$ 
are precisely the index sets $\LP_\ib(y)$ for $y \in W^{\leq w}$.}}

The representation theoretic approach \cite{Y1} to the spectra 
$\Spec \UU^-[w]$ relies on a family of surjective 
$\Tset^r$-equivariant antihomomorphisms $\phi_w \colon R^w_0 \to \UU^-[w]$, 
where $R^w_0$ are certain quotients of subalgebras of the quantum groups 
$R_q[G]$. The algebras $R^w_0$ were introduced by Joseph \cite{J} 
as quantizations 
of the coordinate rings of $w$-translates of the open Schubert cell of 
the flag variety of $\g$, see \S \ref{2.3} for details. Via these maps
one can transfer back and forward questions on the spectra of $\UU^-[w]$ 
to questions on the spectra of quantum function algebras. The latter 
can be approached via representation theoretic methods, building 
on the works of Joseph \cite{J0,J}, Gorelik \cite{Go}, and
Hodges--Levasseur--Toro \cite{HLT}.
This leads to an explicit picture for $\TSpec \UU^-[w]$. 
First, the $\Tset^r$-invariant prime ideals of $\UU^-[w]$
are parametrized by $W^{\leq w}$, and the ideal 
$I_w(y)$ corresponding to $y \in W^{\leq w}$ is explicitly 
given in terms of Demazure modules using the maps $\phi_w$, 
see \thref{Class} for a precise statement. Second, each of the 
strata $\Spec_{I_w(y)} \UU^-[w]$ consists of ideals
constructed by contractions from localizations of $\UU^-[w]/I_w(y)$
by explicit small multiplicative sets of normal elements.

Each of the above two methods has many advantages over the 
other. Using the representation theoretic approach,
it was proved that all ideals $I_w(y)$ are 
polynormal, it was established that $\UU^-[w]$ 
are catenary and satisfy Tauvel's height formula,
the containment problem for $\TSpec \UU^-[w] 
= \{ I_w(y) \mid y \in W^{\leq w} \}$ was solved,
and theorems for separation of variables for $\UU^-[w]$ 
were established (see \cite{Y1,Y5,Y6}). In the special 
case of the algebras of quantum matrices, catenarity and ideal 
containment was proved earlier \cite{Ca2,La} within the framework of 
the ring theoretic approach (though with more complicated 
arguments), but there was no progress on polynormality 
or proofs of these results for more general $\UU^-[w]$ algebras.
On the other hand using the ring theoretic approach, 
it was proved that for all $\Tset^r$-primes $J_D$ 
of $\UU^-[w]$ the factor $\UU^-[w]/J_D$ always has a localization 
that is a quantum torus, its center 
(which is closely related to the structure of the stratum 
$\Spec_{J_D} \UU^-[w]$) was described, and in the case of 
quantum matrices $\Tset^r$-primes were related to total positivity
(see \cite{Ca,BCL,GLL}).
 
Our first result resolves Question 5.3.3 of Cauchon and M\'eriaux
\cite{MC} and unifies the 
two approaches to $\TSpec \UU^-[w]$:
 
\bth{1} Let $\KK$ be an arbitrary base field, 
$q \in \KK^*$ not a root of unity, $\g$ a simple 
Lie algebra, $w$ an element of the Weyl group of $\g$, 
and $\ib$ a reduced word for $w$. Consider the 
presentation of the quantum Schubert cell algebra $\UU^-[w]$
as an iterated Ore extension corresponding to $\ib$.

Then for all Weyl group elements $y \leq w$ 
the Cauchon diagram of the $\Tset^r$-prime ideal 
$I_w(y)$ of $\UU^-[w]$ (from the representation theoretic
approach from \thref{Class} (i)) is equal to 
$\LP_\ib(y)$, the index set of the left positive subword 
of $\ib$ whose total product is $y$.

Thus the $\Tset^r$-prime ideals of $\UU^-[w]$ from the representation 
theoretic approach are related to the ideals $J_D$ from the ring theoretic 
approach via
\[
I_w(y) = J_{\LP_{\ib}(y)}.
\]
\eth 
 
Furthermore, we prove a theorem that explicitly describes the behavior of the 
representation theoretic ideals $I_w(y)$ of $\UU^-[w]$ 
in each stage of the Cauchon deleting derivation procedure. 
This appears in \thref{main2-ind} below and will not be stated 
in the introduction since it requires additional background.
 
With the help of \thref{1}, one can now combine the strengths of the two 
approaches to the spectra of the quantum Schubert cell algebras. We expect 
that the combination of the two methods will lead to substantial progress
in the study of the topology of $\Spec \UU^-[w]$. We use \thref{1} and previous results 
of the second author to resolve Question 5.3.2 of Cauchon and M\'eriaux
\cite{MC}, thereby solving the containment problem for the ideals
\[
\{ J_{\LP_\ib (y)} \mid y \in W^{\leq w} \}
\]
of the classification of \cite{MC}.

\bth{2} In the setting of \thref{1}, the map 
\[
W^{\leq w} \to \TSpec \UU^-[w] \; \; 
\mbox{given by} \; \;
y \mt J_{\LP_{\ib} (y)} 
\] 
is an isomorphism of posets with respect to the
(strong) Bruhat order and the inclusion order on ideals.
\eth
 
Finally, \thref{1} also gives a new, independent proof 
of the Cauchon--Meriaux classification \cite{MC} described 
above. (The proof of \thref{1} does not use results from
\cite{MC}.)
 
Let us return to the general case of Cauchon's method of deleting derivations.
It relates the prime ideals of an initial iterated Ore extension $A$ 
to the prime ideals of the final algebra $\ol{A}$, the Cauchon 
quantum affine space algebra associated to $A$.
In order to study these ideals, one needs an explicit
description of $\ol{A}$ as a subalgebra of the ring 
of fractions $\Fract (A)$. We obtain such for all algebras 
$\UU^-[w]$, establishing yet another relationship 
between the two approaches to the structure of the algebras
$\UU^-[w]$. Given a reduced word 
$\ib=(\al_1, \ldots, \al_l)$ for $w$, define a
successor function $\kap \colon [1,l] \to [1,l] \sqcup \{ \infty \}$
by
\[
\kap(j) = \min \{ k \mid k > j, \al_k = \al_j \}, \; 
\mbox{if} \; \exists k > j \; \; \mbox{such that} 
\; \al_k = \al_j, \quad \kap(j) = \infty, \; 
\mbox{otherwise}.
\]
For $j \in [1,l]$ denote by $\De_{\ib, j} \in \UU^-[w]$ the element obtained 
by evaluating the quantum minor corresponding to the fundamental 
weight $\vpi_{\al_j}$ and the Weyl group 
elements $s_{\al_1} \ldots s_{\al_{j-1}}$, $w$ $\in W$ on the 
$R$-matrix $\RR^w$ corresponding to $w$. We refer to \S \ref{2.3} and 
\S \ref{3.1} for details and the description of these elements in the 
framework of the antiisomorphisms $\phi_w \colon R^w_0 \to \UU^-[w]$.

\bth{3} In the setting of \thref{1}, for all Weyl group 
elements $w$ and reduced words $\ib = (\al_1, \ldots, \al_l)$ 
for $w$, the generators $\ol{x}_1, \ldots, \ol{x}_l$ 
of the corresponding Cauchon quantum affine space algebras 
are given by 
\[
\ol{x}_j = 
\begin{cases}
(q_{\al_j}^{-1} - q_{\al_j})^{-1} 
\De_{\ib, \kap(j)}^{-1} \De_{\ib, j} , 
&\mbox{if} \; \; \kap(j) \neq \infty 
\\
(q_{\al_j}^{-1} - q_{\al_j})^{-1} \De_{\ib, j}, 
&\mbox{if} \; \; \kap(j) = \infty
\end{cases}
\]
for the standard powers $q_{\al_j} \in \KK^*$ of $q$, see \S \ref{2.1}.  
\eth 
This theorem establishes a connection between the initial cluster 
for the cluster algebra structure on $\UU^-[w]$ of
Gei\ss--Leclerc--Schr\"oer and Cauchon's method of deleting 
derivations. We will present a deeper study of this in a 
forthcoming publication. \thref{3} is also an important 
ingredient in a very recent proof \cite{Y-c} of the 
second author of the Andruskiewitsch--Dumas conjecture \cite{AD}.

The paper is organized as follows. Section \ref{qalg}
contains background on the quantum Schubert cell algebras 
and the representation theoretic and ring theoretic 
approaches to the study of their spectra. 
Section \ref{q-aff} contains the proof of \thref{3}. 
Theorems \ref{t1} and \ref{t2} are proved in Section \ref{two-appr}, 
where we also establish a theorem describing the behavior 
of the ideals $I_w(y)$ under the iterations of the 
deleting derivation procedure. 

We will use the following notation throughout the paper.
Given a $\KK$-algebra $A$, we will denote its center by $Z(A)$.
For a $\KK$-subspace $V$ of $A$ and $a, b \in A$ 
we will write $a = b \mod V$ if $a-b \in V$.
Set $\Nset:= \{0,1, \ldots\}$ and $\Zset_+ := \{ 1, 2, \ldots \}$.
For $m, n \in \Zset$ set $[m,n] =\{m, \ldots, n \}$ if $m \leq n$ 
and $[m,n] = \emptyset$ otherwise. 
\medskip
\\
\noindent
{\bf Acknowledgements.} We are thankful to Ken Goodearl for comments 
on the first draft of this paper.
J.G. was supported by the LSU VIGRE NSF 
grant DMS-0739382. M.Y. was supported by NSF grants DMS-1001632 and DMS-1303038. 
\sectionnew{Quantum Schubert cells}
\lb{qalg}
\subsection{Quantized universal enveloping algebras}
\label{2.1}
We will mostly follow the notation of Jantzen's book \cite{Ja}. 
Let $\g$ be a complex simple Lie algebra with root system 
$\Phi$ and Weyl group $W$. Choose a basis $\Pi$ of $\Phi$. Let 
$\lcor.,. \rcor$ be the invariant bilinear form on $\Rset \Pi$
normalized by $\lcor \al, \al \rcor = 2$ for short roots 
$\al \in \Phi$. For $\al \in \Phi$ denote by $\al\spcheck$ and 
$s_\al \in W$ the corresponding coroot and reflection. 
Let $\{ \vpi_\al \mid \al \in \Pi \}$ 
be the fundamental weights of $\g$.  
Denote the root lattice of $\g$ by $\QQ= \Zset \Phi$ and set 
$\QQ^+= \Nset \Phi$. Let $\PP$ be the weight lattice of $\g$
and $\PP^+ = \Nset \{ \vpi_\al \mid \al \in \Pi \}$
be the set of dominant integral weights of $\g$. For a subset
$I \subseteq \Pi$ set $\QQ_I = \Zset I$.
Recall the standard partial order on $\PP$: 
for $\nu_1, \nu_2 \in \PP$ set $\nu_1 \geq \nu_2$ 
if $\nu_2 = \nu_1 - \ga$ for some $\ga \in \QQ^+$. 
Let $\nu_1 > \nu_2$ if $\nu_1 \geq \nu_2$ and $\nu_1 \neq \nu_2$.

Throughout the paper $\KK$ will denote a base field
(of arbitrary characteristic) and $q \in \KK^*$ 
will denote an element which is not 
a root of unity. Denote by $\UU_q(\g)$ the 
quantized universal enveloping algebra of $\g$ 
over $\KK$ with deformation parameter $q$. 
It has generators $K_\al^{\pm 1}, E_\al, F_\al$, 
$\al \in \Pi$ and relations \cite[\S 4.3]{Ja}. The algebra
$\UU_q(\g)$ has a unique Hopf algebra structure 
with comultiplication, antipode, and counit satisfying
\[
\De(K_\al)   = K_\al \o K_\al, \;
\De(E_\al)   = E_\al \o 1 + K_\al \o E_\al, \;
\De(F_\al)   = F_\al \o K_\al^{-1} + 1 \o F_\al
\]
and
\[
S(K_\al) = K^{-1}_\al, \;
S(E_\al)= - K^{-1}_\al E_\al, \; 
S(F_\al)= - F_\al K_\al, \; 
\ep(K_\al)=1, \; \ep(E_\al)= \ep(F_\al)=0.
\]
The subalgebras of $\UU_q(\g)$ generated by 
$\{E_\al\mid \al \in \Pi \}$, 
$\{F_\al\mid \al \in \Pi \}$, 
and 
$\{K^{ \pm 1}_\al\mid \al \in \Pi \}$
will be denoted by $\UU^+$, $\UU^-$, and 
$\UU^0$ respectively.

Denote by $\leq$ the (strong) Bruhat order on $W$
and by $\ell \colon W \to \Nset$ the standard length function.
For $w \in W$ set $W^{\leq w}= \{ y \in W \mid y \leq w \}$.
Let $\B_\g$ be the braid group of $\g$ and $\{T_\al \mid \al \in \Pi \}$
be its standard generating set. We will use Lusztig's action 
of $\B_\g$ on $\UU_q(\g)$ by algebra automorphisms 
in the version given in \cite[\S 8.14]{Ja} by eqs. 
8.14 (2), (3), (7), and (8).  

We will use the following notation for $q$-integers and factorials:
\[
[n]_q := \frac{q^n - q^{-n}}{q- q^{-1}}, \; 
[n]_q ! := [1]_q \ldots [n]_q, \; 
n \in \Nset.
\]
For $\al \in \Pi$, denote $[n]_\al := [n]_{q_\al}$
and $[n]_\al! := [n]_{q_\al}!$,
where $q_\al := q^{\lcor \al, \al \rcor/2}$.
\subsection{Weyl group combinatorics and quantum Schubert cell algebras}
\label{2.2}
Fix $w \in W$. A word $\ib = (\al_1, \ldots, \al_l)$ in the alphabet 
$\Pi$ is called a reduced word for $w$ if $s_{\al_1} \ldots s_{\al_l}$ 
is a reduced expression of $w$ 
(in particular, $\ell(w) = l$). Given
a reduced word $\ib= (\al_1, \ldots, \al_l)$ for $w$,
denote
\begin{equation}
\label{leq1}
w(\ib)_{\leq j} := s_{\al_1} \ldots s_{\al_j} \; \; 
\mbox{and} \; \; w(\ib)_{>j} := s_{\al_{j+1}} \ldots s_{\al_l}
\;\; \mbox{for} \; \; j \in [0,l].
\end{equation}
Thus $w(\ib)_{\leq 0} = 1$ and $w(\ib)_{\leq l} = w$.
There is a bijection between the set of subwords of $\ib$ and 
the subsets of $[1,l]$, which associates 
to a subword $(\al_{j_1}, \ldots, \al_{j_n})$ of $\ib$ 
its index set  
$\{j_1 < \ldots < j_n \} \subseteq[1,l]$.
Given $D \subseteq [1,l]$, 
for $j \in [1,l]$ set 
$s_j^D = s_{\al_j}$ if $j \in D$, and $s_j^D = 1$ otherwise.
Denote 
\begin{equation}
\label{leq2}
w(\ib)^D_{\leq j} := s^D_1 \ldots s^D_j \; \; 
\mbox{and} \; \; w(\ib)^D_{>j}:= s^D_{j+1} \ldots s^D_{l}
\;\; \mbox{for} \; \; j \in [1,l].
\end{equation}
Let
\[
w(\ib)^D := w(\ib)^D_{\leq l } = s^D_1 \ldots, s^D_l.
\]
Following \cite{MR} we call a subword of $\ib$ {\em{(right) positive}} 
if its index set $D \subseteq[1,l]$ has the property that 
\[
w(\ib)^D_{\leq j} 
s_{\al_{j+1}} > w(\ib)^D_{\leq j } \; \; \mbox{for all} \; \; 
j \in [1,l-1].   
\]
A subword of $\ib$ will be called {\em{left positive}} if  
its index set $D \subseteq [1,l]$ has the property
that
\begin{equation}
\label{l-positive}
s_{\al_j} w(\ib)^D_{ > j} > w(\ib)^D_{ > j }\; \; \mbox{for all} \; \; 
j \in [1,l-1].
\end{equation}
Some authors refer to the left positive subwords of $\ib$ 
as Cauchon diagrams associated to 
$\ib$. However, we will use the term Cauchon diagrams for the 
general Cauchon procedure of deleting derivations in iterated 
Ore extensions (see \S \ref{2.4}), and using the same term for different 
notions will easily lead to confusions. 

The map $(\al_{j_1}, \ldots, \al_{j_n} ) \mt ( \al_{j_n}, \ldots, \al_{j_1})$ 
establishes a bijection between the left positive subwords of $\ib$ and 
the right positive subwords of the reduced word $(\al_l, \ldots, \al_1)$ of $w^{-1}$.
Since the map $y \mt y^{-1}$ is a bijection between $W^{\leq w}$ and $W^{\leq w^{-1}}$,
Lemma 3.5 of Marsh--Rietsch \cite{MR} gives that for each $y \in W^{\leq w}$ 
there exists a unique left positive subword of $\ib$ 
such that its index set $D \subseteq [1,l]$ satisfies 
$w(\ib)^D = y$. Denote this index set $D$ by $\LP_{\ib}(y)$.

The support of $w \in W$ is defined by 
\begin{equation}
\label{supp}
\SS(w) := \{ \al \in \Pi \mid s_\al \leq w \}.
\end{equation}
Its complement is given by 
\begin{equation}
\label{compl}
\Pi \backslash \SS(w) 
= \{ \al \in \Pi \mid w \vpi_\al = \vpi_\al \},
\end{equation}
see \cite[Lemma 3.2 and eq. (3.2)]{Y4}.

The quantum Schubert cell algebras $\UU^\pm[w]$, $w \in W$
were defined by De Concini, Kac, and Procesi 
\cite{DKP}, and Lusztig \cite[\S 40.2]{L} as follows. 
Given a reduced word $\ib = (\al_1, \ldots, \al_l)$ 
for $w$, define the roots
\begin{equation}
\label{beta}
\beta_j := w(\ib)_{\leq (j-1)} \al_j, \; \; j \in [1, l]
\end{equation}
and the Lusztig root vectors
\begin{equation}
E_{\be_j} :=
T_{\al_1} \ldots T_{\al_{j-1}} 
(E_{\al_j}), \; \;  
F_{\be_j} := 
T_{\al_1} \ldots T_{\al_{j-1}} 
(F_{\al_j}), \; \; 
j \in [1, l],
\label{rootv}
\end{equation}
see \cite[\S 39.3]{L}.
By \cite[Proposition 2.2]{DKP} and 
\cite[Proposition 40.2.1]{L} the subalgebras 
$\UU^\pm[w]$ of $\UU^\pm$ generated by 
$E_{\beta_j}$, $j \in [1,l]$ and $F_{\beta_j}$, $j \in [1,l]$
do not depend on the choice of a 
reduced word $\ib$ for $w$ and have the PBW bases
\begin{equation}
\label{PBW}
\{(E_{\be_l})^{n_l} \ldots (E_{\be_1})^{n_1} \mid 
n_1, \ldots, n_l \in \Nset
\} \; \; \mbox{and} \; \; 
\{(F_{\be_l})^{n_l} \ldots (F_{\be_1})^{n_1} \mid 
n_1, \ldots, n_l \in \Nset \},
\end{equation}
respectively. 

The algebra $\UU_q(\g)$ is $\QQ$-graded 
by $\deg K_\al = 0$, $\deg E_\al = \al$, $\deg F_\al = - \al$, 
$\forall \al \in \Pi$. This induces a $\QQ$-grading on
$\UU^\pm[w]$. The corresponding graded components will 
be denoted by $(\UU_q(\g))_\ga$ and $(\UU^\pm[w])_\ga$. One has
\begin{equation}
\label{suppU}
\Zset \{ \ga \in \QQ \mid (\UU^\pm[w])_\ga \neq 0 \} = 
\QQ_{\SS(w)}, 
\end{equation}
see e.g. \cite[eq. (2.44) and Lemma 3.2 (ii)]{Y4}.

Recall that there is a unique algebra automorphism $\omega$ 
of $\UU_q(\g)$ such that 
\[
\om(E_\al) = F_\al, \; 
\om(F_\al) = E_\al, \; 
\om(K_\al) = K_{\al}^{-1}, \; \; 
\forall \, \al \in \Pi.
\]
It satisfies $\om( T_\al(u)) = (-1)^{\lcor \al\spcheck, \ga \rcor } 
q^{ - \lcor \al, \ga \rcor }  T_\al( \om (u))$, 
for all $\ga \in \QQ$, $u \in (\UU_q(\g))_\ga$, 
see \cite[eq. 8.14(9)]{Ja}. In other words, 
if $\rho$ is the sum of all fundamental weights of $\g$ and
$\rho\spcheck$ is the sum of all fundamental coweights of $\g$, then 
$\om( T_\al(u)) = (-1)^{ \lcor s_\al(\ga)- \ga, \rho\spcheck \rcor}  
q^{ - \lcor s_\al(\ga) - \ga , \rho \rcor }  T_\al( \om (u))$
for $u \in (\UU_q(\g))_\ga$. Thus
\[
\om( T_y(u)) = (-1)^{ \lcor y(\ga)- \ga, \rho\spcheck \rcor}  
q^{ - \lcor y(\ga) - \ga , \rho \rcor }  T_y( \om (u)), \; \; 
\mbox{for all} \; \;  
y \in W, \ga \in \QQ, u \in (\UU_q(\g))_\ga,
\]
see \cite[eq. 8.18(5)]{Ja} for an equivalent formulation of this fact.
In particular, the 
restrictions of $\om$ induce the isomorphisms
\begin{equation}
\label{om}
\om \colon \UU^+[w] \congto \UU^-[w], 
\quad \om(E_{\be_j}) = 
(-1)^{ \lcor \be_j - \al_j, \rho\spcheck \rcor } 
q^{ - \lcor \be_j - \al_j, \rho \rcor } 
F_{\be_j}, 
\; \; \forall \, j \in [1, \ell(w)].
\end{equation}

To each $\ga \in \QQ$ associate the character
of $\Tset^r = (\KK^*)^{\times \, r}$
\begin{equation}
\label{Tchar}
t \mt 
t^\ga := \prod_{ \al \in  \Pi} t_\al^{\lcor \ga, \vpi_\al \rcor}, \quad 
t =(t_\al)_{\al \in \Pi} \in \Tset^r.
\end{equation}
Define the rational $\Tset^r$-action 
on $\UU_q(\g)$ by algebra automorphisms
\begin{equation}
\label{torus-act}
t \cdot x = t^\ga x, \quad x \in (\UU_q(\g))_\ga.
\end{equation}
It preserves the subalgebras $\UU^\pm[w]$. 
We will denote by $\TSpec  \UU^-[w]$ the space 
of $\Tset^r$-prime ideals of $\UU^-[w]$. 

Fix a reduced word $\ib$ for $w$ and consider the roots \eqref{beta}.  
Eq. \eqref{suppU} implies that for all $j \in [1, \ell(w)]$ 
there exists a unique $t_j = ( t_{j,\al})_{\al \in \Pi} \in \Tset^r$ 
such that 
\begin{equation}
\label{tj}
t_j^{\be_k} = 
q^{\lcor \be_k, \be_j \rcor }, \; \; 
\forall \, k \leq j \quad \mbox{and} \quad
t_{j,\al} = 1, \; \; 
\forall \, \al \in \Pi \backslash \SS( w(\ib)_{\leq j}),
\end{equation}
recall \eqref{Tchar}. The Levendorskii--Soibelman 
straightening law is the following commutation relation 
in $\UU^-[w]$
\begin{multline}
\label{LS}
F_{\be_j} F_{\be_k} - 
q^{ - \lcor \be_k, \be_j \rcor }
F_{\be_k} F_{\be_j}  \\
= \sum_{ {\bf{n}} = (n_{k+1}, \ldots, n_{j-1}) \in \Nset^{\times (j-k-2)} }
p_{\bf{n}} (F_{\be_{j-1}})^{n_{j-1}} \ldots (F_{\be_{k+1}})^{n_{k+1}},
\; \; p_{\bf{n}} \in \KK,
\end{multline}
for all $k < j$, see e.g. \cite[Proposition I.6.10]{BG}. 
The following lemma is a direct consequence
of \eqref{PBW}, \eqref{tj}, and \eqref{LS}.

\ble{Uw} For all base fields $\KK$, $q \in \KK^*$ not a root of unity, 
Weyl group elements $w \in W$ of length $l$, 
reduced words $\ib=(\al_1, \ldots, \al_l)$ 
for $w$, and $j \in [1,l]$ we have:

(i) The subalgebra of $\UU^-[w]$ generated by $F_{\be_1}, \ldots, F_{\be_j}$ 
is equal to $\UU^-[w(\ib)_{\leq j}]$.

(ii) The algebra $\UU^-[w(\ib)_{\leq j}]$ is isomorphic 
to the Ore extension $\UU^-[w(\ib)_{\leq (j-1)}][x_j, \sig_j, \delta_j]$,
where $\sig_j = (t_j \cdot) \in \Aut(\UU^-[w(\ib)_{\leq (j-1)}])$ and
$\de_j$ is a locally nilpotent (left) $\sig_j$-skew derivation
of $\UU^-[w(\ib)_{\leq (j-1)}]$ satisfying $\sig_j \delta_j = 
q_{\al_j}^{-2} \delta_j \sig_j$.
This isomorphism is given by the identity map on 
$\UU^-[w(\ib)_{\leq (j-1)}]$ and $F_{\be_j} \mt x_j$.
Furthermore, 
$\UU^-[w(\ib)_{\leq 0}]= \UU^-[1] \cong \KK$,
$\sig_1=\id$, and $\delta_1 =0$.

(iii) The eigenvalues $t_j \cdot F_{\be_j} = 
q_{\al_j}^{-2} F_{\be_j}$ are not roots of unity.
\ele

The $\sig_j$-skew derivation $\delta_j$ of $\UU^-[w(\ib)_{\leq (j-1)}]$ 
in part (ii) of the lemma is explicitly given by 
\begin{equation}
\label{deltaj}
\delta_j(x) := F_{\be_j} x - q^{\lcor \be_j, \ga \rcor} x F_{\be_j}, \; \; 
\mbox{for} \; \; 
x \in (\UU^-[w(\ib)_{\leq (j-1)}])_\ga, \ga \in \QQ
\end{equation}
and is computed using \eqref{LS}. 

The isomorphisms from part (ii) give rise to the Ore extension
presentations
\[
\UU^-[w(\ib)_{\leq j}] = \UU^-[w(\ib)_{\leq (j-1)}][F_{\be_j}, \sig_j, \delta_j], 
\; \; 1 \leq j \leq l.
\]
When those are iterated, for each reduced word 
$\ib$ for $w$, one obtains a presentation of $\UU^-[w]$ as an 
iterated Ore extension
\begin{equation}
\label{Uwiter}
\UU^-[w] = \KK [F_{\be_1}] [F_{\be_2}; \sig_2, \delta_2] \ldots 
[F_{\be_l}; \sig_l, \delta_l ].
\end{equation}
\subsection{The prime spectrum of $\UU^-[w]$ via Demazure modules}
\label{2.3}
We proceed with the realization of the algebras $\UU^-[w]$ in terms of 
quantum function algebras and the description of the spectra of $\UU^-[w]$
via Demazure modules from \cite{Y1}.

The $q$-weight spaces of a $\UU_q(\g)$-module $V$ 
are defined by 
\[
V_\nu := \{ v \in V \mid K_\al v = q^{ \lcor \nu, \al \rcor} v, \; \; 
\forall \, \al \in \Pi \}, \; \nu \in \PP.
\]
A $\UU_q(\g)$-module is called a type one module if 
$V = \oplus_{\nu \in \PP} V_\nu$.
The category of (left) finite dimensional type one 
$\UU_q(\g)$-modules is semisimple 
(see  \cite[Theorem 5.17]{Ja} and the remark on p. 85 of 
\cite{Ja}). It is closed under taking tensor products and duals
(defined as left modules using the antipode 
of $\UU_q(\g)$). Denote by $V(\la)$ the 
irreducible type one $\UU_q(\g)$-module of highest weight 
$\la \in \PP^+$. Those exhaust all irreducible finite dimensional 
type one modules, see \cite[Theorem 5.10]{Ja}. 

For algebraically closed fields $\KK$ of characteristic 0, 
we will denote by $G$ the connected, simply connected algebraic 
group with Lie algebra $\g$. For all base fields $\KK$ and 
deformation parameters $q \in \KK^*$ that are not roots of unity, 
the quantum group $R_q[G]$ is defined as the Hopf subalgebra 
of the restricted dual $(\UU_q(\g))^\circ$,  
spanned by the matrix coefficients of the modules $V(\la)$, 
$\la \in \PP^+$. The latter are given by 
\begin{equation} 
\label{c-notation}
c^\la_{\xi, v} \in (\UU_q(\g))^\circ,\quad 
c^\la_{\xi, v}(x) := \xi ( x v ), \quad v \in V(\la), 
\xi \in V(\la)^*, x \in \UU_q(\g).
\end{equation}
Because we work with arbitrary base fields, in the notation $R_q[G]$, 
$G$ is just a symbol. 

For each $\la \in \PP^+$, fix a highest weight vector $v_\la$ 
of $V(\la)$. Set for brevity
\[
c_{\xi}^\la := c^\la_{\xi, v_\la}, \; \; \la \in \PP^+, 
\xi \in V(\la)^*.
\]
Define the subalgebra
\[
R^+ := \Span \{ c^\la_\xi \mid \la \in \PP^+, \xi \in V(\la)^* \}
\]
of $R_q[G]$. 

The braid group $\B_\g$ acts on 
the finite dimensional type one 
$\UU_q(\g)$-modules $V$ by 
\begin{equation}
\label{braid}
T_\al(v) := \sum_{l, m, n}
\frac{(-1)^m q_\al^{m-ln}}{[l]_\al! [m]_\al! [n]_\al!} 
 E_\al^{l} F_\al^{m} E_\al^{n} v, 
\quad v \in V_\mu, \mu \in \PP,
\end{equation}
where the sum is over $l,m, n \in \Nset$ such that 
$-l+m-n = \lcor \mu, \al\spcheck \rcor$, cf. 
\cite[\S 8.6]{Ja} and \cite[\S 5.2]{L}.
This action and the $\B_\g$-action on 
$\UU_q(\g)$ are compatible by 
\begin{equation}
\label{compatible}
T_w ( x . v ) := (T_w x) . (T_w v), \; \; \forall \, 
w \in W, x \in \UU_q(\g), v \in V(\la), \la \in \PP^+, 
\end{equation}
see \cite[eq. 8.14 (1)]{Ja}.
Moreover, $T_w(V(\la)_\mu) = V(\la)_{w \mu}$, 
$\forall w \in W$, $\la \in \PP^+$, $\mu \in \PP$.
In particular $\dim V(\la)_{w \la}=1$, 
$\forall w \in W$, $\la \in \PP^+$.

For $\al \in \Pi$ denote by $\UU^\al$ the subalgebra 
of $\UU_q(\g)$ generated by $E_\al$, $F_\al$, and $K_\al^{\pm 1}$:
\begin{equation}
\label{Ual}
\UU^\al = \KK \lcor E_\al, F_\al, K_\al^{\pm 1} \rcor.
\end{equation}
It is canonically isomorphic to $\UU_{q_\al}(\sl_2)$. We will 
later need the following formulas for the irreducible type 
one finite dimensional $\UU^\al$-modules. For all $m, N \in \Nset$, 
$m \leq N$ we have
\begin{equation}
\label{sl2-braid}
T_\al v_{ N \vpi_\al} = \frac{(-q_\al)^N}{[N]_\al !} F_\al^N v_{N \vpi_\al}, \; \;
T_\al^{-1} v_{ N \vpi_\al} = \frac{1}{[N]_\al !} F_\al^N v_{N \vpi_\al},
\end{equation}
and 
\begin{equation}
\label{sl2-EF}
E^m_\al F^m_\al v_{N \vpi_\al} 
= \frac{[m]_\al! [N]_\al!}{[N-m]_\al !} v_{N \vpi_\al},
\end{equation}
by \cite[eqs. 8.6 (6), (7), and Lemma 1.7]{Ja}.

For $\la \in \PP^+$ and $w \in W$ let
$\xi_{w, \la} \in (V(\la)^*)_{- w\la}$ 
be the unique vector such that
\begin{equation}
\label{xi-w}
\lcor \xi_{w, \la}, T^{-1}_{w^{-1}} v_\la \rcor =1. 
\end{equation}
For $y, w \in W$ and $\la \in \PP^+$ define the quantum minors
\begin{equation}
\label{e}
e^\la_{y,w} := c^\la_{\xi_{y,\la}, T_{w^{-1}}^{-1} v_\la} \in R_q[G] \; \; 
\mbox{and} \; \; 
e^\la_w := e^\la_{1,w} = c^\la_{\xi_{w,\la}} \in R^+.
\end{equation}
Using the second equality in \eqref{sl2-braid} one easily shows that they coincide with 
the quantum minors of Berenstein and Zelevinsky from \cite[Eq. (9.10)]{BZ}. 
If one works with $T_w$ instead of $T_{w^{-1}}^{-1}$, then 
additional scalars arise from the first equality in \eqref{sl2-braid}.
This is why we use the latter throughout the paper.

We have 
\begin{equation}
\label{mult}
e^{\la_1}_w e^{\la_2}_w = e^{\la_1 + \la_2 }_w 
= e^{\la_2}_w e^{\la_1}_w,  \; \;  \forall \, \la_1, \la_2 \in \PP^+, 
w \in W,
\end{equation}
which is proved analogously to \cite[eq. (2.18)]{Y4}
using one more time the second equality in \eqref{sl2-braid}. 
Joseph proved that the multiplicative sets 
$E_w = \{ e_w^\la \mid \la\in \PP^+ \} \subset R^+$ are Ore sets, see
\cite[Lemma 9.1.10]{J}. Joseph's proof works for all base fields 
$\KK$, $q \in \KK^*$ not a root of unity, see \cite[\S 2.2]{Y6}.
Define the quotient algebras
\[
R^w := R^+ [E_w^{-1}], \quad w \in W
\]
and their subalgebras
\begin{equation}
\label{R0w}
R^w_0 := \{ c^\la_\xi (e^\la_w)^{-1} \mid \la \in \PP^+, \xi \in V(\la)^* \},
\end{equation}
introduced by Joseph \cite[\S 10.4.8]{J}.
One does not need to take span in the right hand side 
of the above formula, cf. \cite[\S 10.4.8]{J} or \cite[eq. (2.18)]{Y5}.
The algebra $R^w_0$ is $\QQ$-graded by 
\[
(R_0^w)_\ga := \{  c^\la_\xi (e^\la_w)^{-1} \mid \la \in \PP^+, 
\xi \in (V(\la)^*)_{\ga + w(\la)} \}, \; \; 
\ga \in \QQ.
\]
For $\mu = \la_1 - \la_2 \in \PP$, $\la_1, \la_2 \in \PP^+$, set
\begin{equation}
\label{emu}
e^\mu_w := e^{\la_1}_w (e^{\la_2}_w)^{-1} \in R^w_0.
\end{equation}
It follows from \eqref{mult} that this does not depend on the choice 
of $\la_1, \la_2$ and that 
$e^{\mu_1}_w e^{\mu_2}_w = e^{\mu_1 + \mu_2}_w$ 
for all $\mu_1, \mu_2 \in \PP$.

The $\UU^\pm \UU^0$-submodules $\UU^\pm V(\la)_{y \la} = \UU^\pm T_y v_\la$
of $V(\la)$, where $y \in W$, are called Demazure modules. They give rise to the 
quantum Schubert cell ideals of $R^+$
\[
Q(y)^\pm := \Span \{ c^\la_\xi \mid \la \in \PP^+, \xi \in V(\la)^*, \,
\xi \perp \UU^\pm T_y v_\la \}, \; \; y \in W. 
\]
Their counterparts in the algebras $R^0_w$ are the ideals 
\begin{equation}
\label{Qyw}
Q(y)^\pm_w := \{ c^\la_\xi e^{-\la}_w \mid \la \in \PP^+, 
\xi \in V(\la)^*, \,
\xi \perp \UU^\pm T_y v_\la \} = Q(y)^\pm E_w^{-1}
\cap R_0^w. 
\end{equation}
Analogously to \eqref{R0w} one does not need 
to take a span in \eqref{Qyw}, see \cite{Go,Y1}.
For $\ga \in \QQ^+ \backslash \{ 0 \}$ denote
$m_w(\ga) = \dim (\UU^+[w])_\ga= \dim (\UU^-[w])_{-\ga}$. Let
$\{u_{\ga, i} \}_{i=1}^{m_w(\ga)}$ and 
$\{u_{-\ga, i} \}_{i=1}^{m_w(\ga)}$ be 
dual bases of $(\UU^+[w])_\ga$ and 
$(\UU^-[w])_{-\ga}$ with respect to the Rosso--Tanisaki form,
see \cite[Ch. 6]{Ja}. The quantum $R$ matrix corresponding 
to $w$ is given by
\begin{equation}
\label{RR}
\RR^w := 1 \otimes 1 + \sum_{\ga \in \QQ^+, \ga \neq 0} \sum_{i=1}^{m_w(\ga)} 
u_{\ga, i} \otimes u_{- \ga, i} \in \UU^+ \wh{\otimes} \UU^-,
\end{equation}
where $\UU^+ \wh{\otimes} \UU^-$ is the completion of $\UU^+ \otimes \UU^-$ 
with respect to the descending filtration \cite[\S 4.1.1]{L}.
Finally, we recall that there is a unique graded algebra 
antiautomorphism $\tau$ of $\UU_q(\g)$ defined by 
\begin{equation}
\label{tau}
\tau(E_\al) = E_\al,
\,
\tau(F_\al) = F_\al, 
\, 
\tau(K_\al) = K_\al^{-1}, \; \; 
\al \in \Pi,
\end{equation}
see \cite[Lemma 4.6(b)]{Ja}. It satisfies
\begin{equation}
\label{tau-ident}
\tau (T_w x) = T_{w^{-1}}^{-1} ( \tau (x)), \; \; 
\forall \, x \in \UU_q(\g), w \in W,
\end{equation}
see \cite[ eq. 8.18(6)]{Ja}. 

The next theorem summarizes the representation theoretic 
approach to $\Spec \UU^-[w]$ via quantum function algebras
and Demazure modules.

\bth{Class} 
For all base fields $\KK$, $q \in \KK^*$ not a root of unity, 
simple Lie algebras $\g$, and Weyl group elements $w \in W$, 
the following hold:

(i) The maps 
\begin{equation}
\label{phi}
\phi_w \colon 
R^w_0 \to \UU^-[w], \quad
\phi_w \big( c^\la_\xi e^{-\la}_w \big) 
:= \big( c^\la_{\xi, T^{-1}_{w^{-1}} v_\la} \otimes \id \big) (\tau \otimes \id) 
\RR^w, \; \; 
\la \in \PP^+, \xi \in V(\la)^*
\end{equation}
are well defined surjective $\QQ$-graded algebra antihomomorphisms
with kernels $\ker \phi_w = Q(w)_w^+$. 

(ii) For $y \in W^{\leq w}$ the ideals
\[
I_w(y) = \phi_w( Q(w)_w^+ +  Q(y)_w^-) = \phi_w(Q(y)_w^-)
\]
are distinct, $\Tset^r$-invariant, completely prime  
ideals of $\UU^-[w]$. All $\Tset^r$-prime ideals of $\UU^-[w]$
are of this form.

(iii) The map $y \in W^{\leq w} \mt I_w(y) \in \TSpec \UU^-[w]$ 
is an isomorphism of posets with respect to the Bruhat order 
on $W^{\leq w}$ and the inclusion order on $\TSpec \UU^-[w]$.
\eth

Part (i) is \cite[Theorem 2.6]{Y4}. It was first proved in \cite{Y1} 
for another version of the Hopf algebra $\UU_q(\g)$ equipped with 
the opposite coproduct, a different braid group action and 
Lusztig's root vectors.
Theorem 2.6 in \cite{Y4} used $T_w$ in place 
of $T^{-1}_{w^{-1}}$ in eqs. \eqref{xi-w} and \eqref{phi}. The two 
formulations are equivalent since $\dim V(\la)_{w \la} = 1$ 
and $T_w(V(\la)_\mu) = V(\la)_{w \mu}$ for all $w \in W$,
$\la \in \PP^+$, $\mu \in \PP$. Parts (ii)--(iii) of \thref{Class}
are proved in \cite[Theorem 3.1 (a)]{Y6} relying on results 
of Gorelik \cite{Go} and Joseph \cite{J0}. These statements 
were earlier proved in \cite[Theorem 1.1 (a)-(b)]{Y1} under slightly stronger 
conditions on $\KK$ and $q$. 

Recall \eqref{e}. The elements 
\[
b_{y,w}^\la := \phi_w( e^\la_y e^{-\la}_w )
= (e^\la_{y,w} \tau \otimes \id) \RR^w, 
\quad \la \in \PP^+
\]
are nonzero normal elements of $\UU^-[w]/I_w[y]$:
\begin{equation}
\label{commute}
b_{y,w}^\la x
= 
q^{ - \lcor (w+y)\la, \ga \rcor }
x b_{y,w}^\la, \; \; 
\forall \, 
\la \in \PP^+,
\ga \in \QQ_{\SS(w)}, 
x \in (\UU^-[w]/I_w(y))_\ga,
\end{equation}
by \cite[Theorem 3.1(b) and eq. (3.1)]{Y5}. 
Here and below we denote by the 
same symbols the images of elements of $\UU^-[w]$ 
and $R_q[G]$ in their factors, which
is a standard notational convention. The $R$-matrix
commutation relations in $R^+$ (see e.g. \cite[Theorem I.8.15]{BG})
and eq. \eqref{mult} imply that for all 
$\la_1, \la_2 \in \PP^+$, $b_{y,w}^{\la_1}b_{y, w}^{\la_2} = 
q^{- \lcor \la_1, \la_2 - y^{-1} w \la_2 \rcor}
b_{y,w}^{\la_1+ \la_2}$. 
Thus 
\[
B_{y,w} := \KK^*\{ b_{y,w}^\la \mid \la \in \PP^+ \}
\]
is a multiplicative subset of $\UU^-[w]/I_w(y)$ consisting of
normal elements. The quotient ring 
$R_{y,w}:= (\UU^-[w]/I_w(y))[B_{y,w}^{-1}]$ is $\Tset^r$-simple. 
Its center is a Laurent polynomial ring of dimension $\dim \ker(w+y)$. 
The prime spectrum of $\UU^-[w]$ is partitioned into
\[
\Spec \UU^-[w] = \bigsqcup_{y \in W^{\leq w}} \Spec_{I_w(y)} \UU^-[w], 
\]
where
\[
\Spec_{I_w(y)} \UU^-[w] := \{ J \in \Spec \UU^-[w] \mid
J \supseteq I_w(y) \; \; \mbox{and} \; \; 
J \cap B_{y,w} = \emptyset \}.
\]
Moreover, extension and contraction establishes the 
homeomorphisms: 
\[
\Spec Z(R_{y,w}) \congto \Spec R_{y,w} \congto 
\Spec_{I_w(y)} \UU^-[w]
\]
and the centers $Z(R_{y,w})$ are Laurent polynomial rings.
We refer to \cite[Theorem 3.1 and Proposition 4.1]{Y5} for details
and proofs of the above statements. The dimensions of the 
Laurent polynomial rings $Z(R_{y,w})$ were explicitly determined 
in \cite{BCL,Y6}.
The above results fit to the general 
framework of Goodearl and Letzter \cite{GL} for reconstruction of 
the spectra of algebras from their torus invariant prime spectra.
Compared to \cite{GL}, the above framework for $\Spec \UU^-[w]$ 
is much more explicit. It deals with explicit $\Tset^r$-prime ideals 
and localizations by small sets of normal elements.

The antihomomorphisms $\phi_w \colon R_0^w \to \UU^-[w]$
are explicitly given by 
\begin{multline}
\label{phi1} 
\phi_w(c_\xi^\la e_w^{-\la} ) = \sum_{m_1, \ldots, m_l \in \Nset}
\left( \prod_{j=1}^l 
\frac{ (q_{\al_j}^{-1} - q_{\al_j})^{m_j}}
{q_{\al_j}^{m_j (m_j-1)/2} [m_j]_{\al_j}! } \right) 
\\
\times 
\lcor \xi, (\tau E_{\be_1})^{m_1} \ldots (\tau E_{\be_l})^{m_l} T^{-1}_{w^{-1}} v_\la \rcor 
F_{\be_l}^{m_l} \ldots F_{\be_1}^{m_1},
\end{multline}
for all $\la \in \PP^+$, $\xi \in V(\la)^*$. This 
follows from \eqref{phi} and the standard formula 
\cite[eqs. 8.30 (1) and (2)]{Ja} for
the inner product of the  pairs of monomials 
\eqref{PBW} with respect 
to the the Rosso--Tanisaki form. 
\subsection{Cauchon's method of deleting derivations}
\label{2.4}
We continue by outlining Cauchon's ring theoretic approach to the study 
of $\Spec \UU^-[w]$ via the method of deleting derivations.
We follow \cite{Ca,MC} and the review in 
\cite[Section 2]{BL}.

Fix an iterated Ore extension 
\begin{equation} 
\label{itOre}
A := \KK[x_1][x_2; \sig_2, \delta_2] \ldots [x_l; \sig_l, \delta_l].
\end{equation}
In particular, for $j \in [2,l]$, $\sig_j$ is an automorphism 
and $\delta_j$ is a (left) $\sig_j$-skew derivation of 
the $(j-1)$-st algebra $A_{j-1}:=\KK[x_1][x_2; \sig_2, \delta_2] 
\ldots [x_{j-1}; \sig_{j-1}, \delta_{j-1}]$
in the above chain. 

\bde{CGL} An iterated Ore extension $A$ as in \eqref{itOre} 
is called a Cauchon--Goodearl--Letzter (CGL) extension if it is 
equipped with a rational action of the torus $\Tset^r= (\KK^*)^{\times r}$, $r \in \Zset_+$
by algebra automorphisms satisfying the following conditions:

(i) The elements $x_1,\dots,x_l$ are $\Tset^r$-eigenvectors.

(ii) For every $j \in [2,l]$, $\de_j$ is a locally nilpotent 
$\sig_j$-derivation of $A_{j-1}$. 

(iii)  For every $j \in [1,l]$, there exists $t_j \in \Tset^r$ such that 
$\sig_j = (h_j \cdot)$ as elements of $\Aut(A_{j-1})$
and the $t_j$-eigenvalue of $x_j$, to be denoted by $q_j$, is not a root 
of unity.
\ede 

One easily deduces that for all CGL extensions, 
$\sig_j \delta_j =q_j \delta_j \sig_j$, $\forall j \in [2,l]$.
For $1 \leq i < j \leq l$ denote the eigenvalues
\[
t_j \cdot x_i = q_{j,i} x_i.
\]
Given a CGL extension $A$ as in \eqref{itOre},
for $j = l+1, l, \ldots, 2$, Cauchon iteratively 
constructed in \cite{Ca} $l$-tuples of nonzero elements 
\[
(x^{(j)}_1, \ldots, x^{(j)}_l)
\] 
and families of subalgebras 
\[
A^{(j)} := \KK \lcor x^{(j)}_1, \ldots, x^{(j)}_l \rcor
\]
of the division ring of fractions $\Fract(A)$ of $A$. 
First, set
\[
(x^{(l+1)}_1, \ldots, x^{(l+1)}_l) := 
(x_1, \ldots, x_l) \; \; \mbox{and} \; \; 
A^{(l+1)}=A.
\]
For $j=l, \ldots, 2$, the $l$-tuple $(x^{(j)}_1, \ldots, x^{(j)}_l)$ 
is determined from
$(x^{(j+1)}_1, \ldots, x^{(j+1)}_l)$ by 
\begin{equation}
\label{new-x}
x^{(j)}_i := 
\begin{cases}
x^{(j+1)}_i, 
& \mbox{if} \; \; i \geq j 
\\
\sum_{m=0}^\infty \frac{(1- q_j)^{-m}}{(m)_{q_j}!} 
\Big[ \delta_j^m \sig^{-m}_j \left(x^{(j+1)}_i \right) \Big]
\left(x^{(j+1)}_j \right)^{-m}, 
& \mbox{if} \; \; i < j.
\end{cases}
\end{equation}
Here $(0)_q=1$, $(m)_q = (1-q^m)/(1-q)$ for $m > 0$, and 
$(m)_q! = (0)_q \ldots (m)_q$ for $m \in \Nset$. 
For $j \in [2,l+1]$, Cauchon constructed an algebra 
isomorphism 
\begin{equation}
\label{Aj-isom}
A^{(j)} \congto
\KK[y_1] \ldots [y_{j-1}; \sig_{j-1}, \delta_{j-1}] 
[y_j; \tau_j] \ldots [y_l; \tau_l],
\end{equation}
where $\tau_k$ denotes the automorphism of 
$\KK[y_1] \ldots [y_{j-1}; \sig_{j-1}, \delta_{j-1}] 
[y_j; \tau_j] \ldots [y_{k-1}; \tau_{k-1}]$ 
such that $\tau_k(y_i) = q_{k,i} y_i$ for all $i \in[1,k-1]$.
This isomorphism is given by 
$x^{(j)}_i \mt y_i$, $i= 1, \ldots, l$. Define 
\[
S_j := \Big\{ \left(x^{(j+1)}_j \right)^m 
\; \Big{|} \; m \in \Nset \Big\}, \; \; j \in [2,l] .
\]
Then $S_j$ is an Ore subset of $A^{(j)}$ and $A^{(j+1)}$.
Cauchon proved that $A^{(j)}[S_j^{-1}] = A^{(j+1)} [S_j^{-1}]$.

Set $q_{i,i}= 1$ for $i \in [1,l]$ and $q_{i,j} = q^{-1}_{j,i}$ 
for $1 \leq i < j \leq l$. 
The quantum affine space algebra $R_{\qb}[\AA^l]$ associated to the 
matrix $\qb := (q_{i,j})_{i,j=1}^l$ is the $\KK$-algebra
with generators $y_1, \ldots, y_l$ and relations 
$y_i y_j = q_{i,j} y_j y_i$, $\forall i, j \in [1,l]$. 
We will call the algebra $A^{(2)}$ obtained at the 
end of the Cauchon deleting derivation 
procedure {\em{the Cauchon quantum affine space algebra associated to}} $A$
and will denote it by $\ol{A}:= A^{(2)}$. Correspondingly, 
the final $l$-tuple of $x$-elements will be denoted by 
$(\ol{x}_1, \ldots, \ol{x}_l) = (x^{(2)}_1, \ldots, x^{(2)}_l)$. 
For $j=2$ eq. \eqref{Aj-isom} gives an isomorphism 
\begin{equation}
\label{final-isom}
\ol{A}= A^{(2)} \congto R_{\qb}[\AA^n], \; \; 
\ol{x}_i = x^{(2)}_i \mt y_i, i \in [1,l]. 
\end{equation}
Furthermore, Cauchon constructed set-theoretic embeddings
\[
\vp_j \colon \Spec A^{(j+1)} \hra \Spec A^{(j)}, \; \; 
j \in [2,l], 
\]
which have certain topological properties 
but are not topological embeddings. They are given by 
\[
\vp_j(J_{j+1})= 
\begin{cases}
J_{j+1} S_j^{-1} \cap A^{(j)}, 
& \mbox{if} \; \; x_j^{(j+1)} \notin J_{j+1} \\
g_j^{-1} \left( J_{j+1}/ \left( x_j^{(j+1)} \right) \right), 
& \mbox{if} \; \; x_j^{(j+1)} \in J_{j+1},  
\end{cases} 
\]
where $J_{j+1} \in \Spec A^{(j+1)}$. Here 
$g_j \colon A^{(j)} \to A^{(j+1)}/(x_j^{(j+1)})$ 
is the homomorphism given by 
$g_j( x_i^{(j)} ) := x_i^{(j+1)} + (x_j^{(j+1)})$, $i \in [1,l]$.
For this construction one needs \cite{Ca} the additional condition 
$x_j^{(j+1)} \notin J_{j+1} \Rightarrow 
J_{j+1} \cap S_{j+1} = \emptyset$.
This condition is satisfied for all $J_{j+1} \in \TSpec A^{(j+1)}$
since by a result of Goodearl and Letzter \cite[Proposition 4.2]{GL} 
all $\Tset^r$-prime ideals of a CGL extension are completely 
prime (recall \eqref{Aj-isom}). 
A CGL extension $A$ as in \deref{CGL} is called torsion free, 
if the subgroup of $\KK^*$ generated by 
all $q_{j,i}$, $1 \leq i < j \leq l$ is torsion free. 
By another result of Goodearl and Letzter \cite[Theorem 2.3]{GL0} 
all prime ideals of a torsion free CGL extension are completely 
prime. Thus the above mentioned condition is satisfied for all 
torsion free CGL extensions $A$ because of \eqref{Aj-isom}.
By \leref{Uw} all algebras $\UU^-[w]$ are 
torsion free CGL extensions when $q \in \KK^*$ 
is not a root of unity.

The composition 
$\vp := \vp_2 \ldots \vp_l \colon \Spec A \hra \Spec \ol{A}$ 
is a set-theoretic embedding, which restricts to an 
embedding $\TSpec A \hra \TSpec \ol{A}$.
Since $\ol{A}$ is 
a quantum affine space algebra, see \eqref{final-isom}, 
the $\Tset^r$-prime ideals of $\ol{A}={A^{(2)}}$ are the 
ideals $K_D := \ol{A} \{ \ol{x}_i \mid i \in D \}$ 
for $D \subseteq [1,l]$. The Cauchon
diagram of $J \in \TSpec A$ 
is the unique set $D \subseteq [1,l]$ such that 
$\vp(J) = K_D$. We will denote the Cauchon diagram of $J$
by $\CD(J)$. If $D \subseteq [1,l]$ is the 
Cauchon diagram of a $\Tset^r$-invariant 
prime ideal of $A$, then this prime ideal 
will be denoted by
\begin{equation}
\label{JD}
J_D:=\vp^{-1}(K_D).
\end{equation}

Let 
\begin{equation}
\label{Aprime}
A' := \KK \lcor x_1, \ldots, x_{l-1} \rcor = 
\KK[x_1][x_2; \sig_2, \delta_2] \ldots [x_{l-1}; \sig_{l-1}, \delta_{l-1}]. 
\end{equation}
So $A = A'[x_l ; \sig_l, \delta_l]$. Set 
\[
A'' = \KK \lcor x_1^{(l)}, \ldots, x_{l-1}^{(l)} \rcor. 
\]
So $A^{(l)} = A''[x_l ; \tau_l]$.
Note that $A'$ and $A''$ are $\Tset^r$-stable subalgebras of $A=A^{(l+1)}$ and $A^{(l)}$,
respectively. They are isomorphic via the following
$\Tset^r$-equivariant algebra isomorphism (recall \eqref{new-x}): 
\begin{equation}
\label{theta}
\theta \colon A' \congto A'', \; \; \theta(a')= 
\sum_{m=0}^\infty \frac{(1- q_l)^{-m}}{(m)_{q_l}!} [ \delta_l^m \sig^{-m}_l(a')] 
x_l^{-m}.
\end{equation}
It satisfies
$\theta( x_i) = x_i^{(l)}$, $i \in [1, l-1]$.
For an ideal $J$ of $A$ denote its leading part consisting 
of the leading terms of the elements of $J$ written as left or right 
polynomials in $x_l$ with coefficients in $A'$:
\begin{align}
\label{hI}
\lt(J) :&= \{ a' \in A' \mid \exists a \in J, \; m \in \Nset \; \; 
\mbox{such that} \; \; a -  a' x_l^m \in A' x_l^{m-1} + \ldots + A' \}
\\
&= \{ a' \in A' \mid \exists a \in J, \; m \in \Nset \; \; 
\mbox{such that} \; \; a -  x_l^m a' \in x_l^{m-1} A' + \ldots + A' \}.
\nn
\end{align}
(The equality holds because $\sig_l$ is locally finite.)

The proof of the following lemma is analogous to \cite[Lemma 4.7]{KL}
and is left to the reader.

\ble{GKdim} Let $x$ be a regular element of the $\KK$-algebra $A$
for which there exist two $\KK$-linear maps $\sig, \delta \colon A \to A$
such that $\sig$ is locally finite, $\delta$ is locally nilpotent, 
$\sig \delta = q \delta \sig$ for some $q \in \KK^*$, and
\[
x a = \sig(a) x + \delta(a), \; \; 
\forall \, a \in A.
\]
Then the set $\Om=\{ 1, x, x^2, \ldots \}$ is an Ore subset 
of $A$ and
\[
\GKdim (A [\Om^{-1}]) = \GKdim A.
\]
\ele

We will need the following facts for a recursive 
computation of Cauchon diagrams and Gelfand--Kirillov 
dimensions of quotients.

\bpr{ind} Assume that $J$ is a $\Tset^r$-prime ideal of a CGL extension 
$A$ given by \eqref{itOre}. 

(i) If $x_l \notin J$, then 
\begin{equation}
\label{ideals}
J S_l^{-1} = \oplus_{m \in \Zset} \theta( \lt(J)) x^m_l, \quad
\vp_l(J)= \oplus_{m \in \Nset} \theta(\lt(J)) x_l^m,
\end{equation}
$\CD(J) = \CD( \lt(J) )$, and 
\begin{equation}
\label{GK1}
\GKdim (A/J) = \GKdim (A'/\lt(J))+1.
\end{equation}

(ii) If $x_l \in J$, then $\vp_l(J) = \theta(J \cap A') + A^{(l)} x_l$,  
$\CD(J) = \CD(J \cap A') \sqcup \{l\}$, and we have the 
$\Tset^r$-equivariant algebra isomorphisms
$A/J \cong A^{(l)}/ \vp_l(J) \cong A'/(J \cap A') \cong A''/(\vp_l(J) \cap A'')$.    
In particular, $\GKdim (A/J) = \GKdim (A'/J \cap A')$.
\epr
Here the Cauchon diagrams $\CD( \lt(J) )$ and $\CD(J \cap A')$ 
are computed with respect to the presentation \eqref{Aprime} 
of $A'$ as a CGL extension. 
\begin{proof} Part (i): By \cite[Lemma 2.2]{LLR}
every $\Tset^r$-invariant ideal $L$ of 
$A S_l^{-1} = A''[x_l^{\pm 1} ; \tau_l]$ has the 
form 
\begin{equation}
\label{Jid}
L = \oplus_{m \in \Zset} L_0 x_l^m \; \; 
\mbox{for some ideal} \; L_0 \;  \mbox{of} \;  A''.
\end{equation}
If $a = \sum_m a_m x_l^m   \in L$, 
then $t_l^k \cdot (x_l^{-k} a x_l^k) 
= \sum_m q_l^{k m} a_m x_l^m  \in L$ for all $k \in \Nset$, 
where $t_l \in \Tset^r$ is the element from \deref{CGL} (iv).
Thus $a_m x_l^m \in L$, $\forall m \in \Zset$, which proves \eqref{Jid}.

We apply this to the ideal $L:= J S_l^{-1}$.
Eq. \eqref{theta} implies that for all $a' \in A'$ and $m \in \Zset$
\[
\theta(a') x^m_l = a' x^m_l + \sum_{k=n}^{m-1} b'_k x^k_l 
\]
for some $n < m$, $b'_k \in A'$. Since every nonzero element of 
$JS_l^{-1}$ has the form $a' x^m_l + \sum_{k =n}^{m-1} a'_k x^k_l$ 
for some $a' \in \lt(I) \backslash \{0\}$, $n <m \in \Zset$, and $a'_k \in A'$
it should also have the form $\theta(a') x^m_l + \sum_{k =n}^{m-1} a''_k x^k_l$ 
for some $a' \in \lt(I) \backslash \{0\}$, $n <m \in \Zset$, and $a''_k \in A''$.
Now the two equalities in \eqref{ideals} follow from \eqref{Jid}. 
The equality $\CD(I) = \CD( \lt(I) )$  
is a consequence of the definition of Cauchon diagrams.
The last statement of part (i) follows from \leref{GKdim} and 
the fact that 
$(A/J)[S_l^{-1}] \cong \theta(A'/\lt(J))[x_l^{\pm}, \tau_l]$. 

Part (ii): The first two statements follow from the definition 
of $\vp_l$. The latter also implies that $g_l$ induces the $\Tset^r$-equivariant 
algebra isomorphism $A^{(l)}/\vp_l(J) \cong A/J$. Since $x_l \in J$ and
$x_l \in \vp_l(J)$ the embeddings $A' \hra A$ and 
$A'' \hra A^{(l)}$ induce the $\Tset^r$-equivariant 
algebra isomorphisms $A'/(J \cap A') \cong A/J$ and 
$A''/(\vp_l(J) \cap A'') \cong A^{(l)}/\vp_l(J)$.
\end{proof}

By \leref{Uw}, the quantum Schubert cell algebras $\UU^-[w]$ are 
torsion free CGL extensions for all base fields $\KK$ and 
$q \in \KK^*$ not a root of unity. 
There is one presentation \eqref{Uwiter} 
of $\UU^-[w]$ as a CGL extension 
for each reduced word $\ib$ for $w$.
Cauchon and M\'eriaux established in \cite{MC} the 
following classification result for their $\Tset^r$-spectra.

\bth{CM} (Cauchon--M\'eriaux, \cite{MC}) For all base fields $\KK$, 
$q \in \KK^*$ not a root of unity, simple Lie algebras $\g$, 
Weyl group elements $w$, and reduced words $\ib$ for $w$,
consider the presentation \eqref{Uwiter} of $\UU^-[w]$ as a 
torsion free CGL extension. In this presentation,
the $\Tset^r$-prime ideals of $\UU^-[w]$ are the ideals 
$J_{\LP_\ib(y)}$ for the elements $y \in W^{\leq w}$
(recall \eqref{JD}),
where $\LP_{\ib}(y) \subseteq [1,l]$ is the index set of 
the left positive subword of $\ib$ whose total product is $y$,
cf. \S \ref{2.2}.
\eth   

In other words the theorem asserts that 
the Cauchon diagrams of the $\Tset^r$-invariant 
prime ideals of $\UU^-[w]$ for the presentation \eqref{Uwiter}
as an iterated Ore extension are precisely the index sets 
of all left positive subwords of $\ib$.
In \cite{MC} \thref{CM} was formulated for the algebras $\UU^+[w]$. 
The two statements are equivalent because of the isomorphism 
\eqref{om}.

We give a second, independent proof of this theorem in Section \ref{two-appr}.
\sectionnew{Cauchon's affine space algebras associated to $\UU^-[w]$} 
\label{q-aff}
\subsection{Statement of main result}
\label{3.1} For each reduced word $\ib$ for a Weyl group element 
$w \in W$ we have a presentation \eqref{Uwiter} of the quantum Schubert cell 
algebra $\UU^-[w]$ as a torsion free CGL extension.
The Cauchon quantum affine space algebra associated to each of the 
algebras $\UU^-[w]$ and a presentation of $\UU^-[w]$ as a CGL extension via a 
reduced words $\ib$ for $w$ is the result of an intricate
iterative procedure. In this section we obtain an explicit description 
of each of these quantum affine space algebras using the 
antiisomorphisms from \thref{Class} (i). This is done in 
\thref{main1}. It expresses each of the generators of the 
Cauchon quantum affine space algebras associated to $\UU^-[w]$ and $\ib$
as a quantum minor or as a fraction of two quantum minors.

Fix a Weyl group element $w \in W$ and a reduced word 
$\ib= (\al_1, \ldots, \al_l)$ for it where $l= \ell(w)$.
Let 
\[
\ol{F}_{\ib, 1}, \ldots, \ol{F}_{\ib, l} 
\]
denote the generators $\ol{x}_1, \ldots, \ol{x}_l$ of the Cauchon 
quantum affine space algebra associated to the presentation \eqref{Uwiter} 
of $\UU^-[w]$ as a CGL extension corresponding to the reduced word 
$\ib$ for $w$, recall \S \ref{2.4}. Define a successor function 
$\kap \colon [1,l] \sqcup \{ \infty \} \to [1,l] \sqcup \{ \infty \}$
associated to $\ib$ as follows.
Let $j \in [1,l]$. If there exists $k > j$ such that $\al_k = \al_j$, 
then we let 
\begin{equation}
\label{s}
\kap(j) = \min \{ k \mid k > j, \al_k = \al_j \}.
\end{equation}
Otherwise, we let $\kap(j) = \infty$. Set $\kap(\infty) = \infty$. Let 
\begin{equation}
\label{m}
O(j) = \max \{ n \in \Nset \mid \kap^n(j) \neq \infty \}, 
\end{equation}
where as usual $\kap^0 := \id$. Define the quantum minors 
\begin{align}
\label{minors}
\De_{\ib, j} :&= b^{\vpi_{\al_j}}_{ w(\ib)_{\leq (j -1)}, w}
=
\phi_w \left( 
e^{\vpi_{\al_j}}_{w(\ib)_{\leq (j -1)}}
e^{-\vpi_{\al_j}}_w 
\right)
\\
&= \left( e^{\vpi_{\al_j}}_{ w(\ib)_{\leq (j -1)}, w} \tau \otimes \id \right) 
\RR^w \in \UU^-[w], 
\; \; j \in [1,\ell(w)],   
\nn
\end{align}
recall \eqref{e}, \eqref{RR}, \eqref{tau}, and \thref{Class} (i).

\bth{main1} Assume that $\KK$ is an arbitrary base field, 
$q \in \KK^*$ is not a root of unity, $\g$ is a simple 
Lie algebra, $w \in W$ is a Weyl group element, and 
$\ib$ is a reduced word for $w$. Then the generators 
$\ol{F}_{\ib,1}, \ldots, \ol{F}_{\ib, \ell(w)}$ 
of the Cauchon quantum affine space algebra associated 
to the presentation \eqref{Uwiter} of $\UU^-[w]$ 
as a CGL extension corresponding to $\ib$ 
are given by
\[
\ol{F}_{\ib, j} = (q_{\al_j}^{-1} - q_{\al_j})^{-1} 
\De_{\ib, \kap(j)}^{-1} \De_{\ib, j} , \; \; 
\mbox{if} \; \; \kap(j) \neq \infty 
\]
and 
\[
\ol{F}_{\ib, j} = (q_{\al_j}^{-1} - q_{\al_j})^{-1} \De_{\ib, j}, \; \; 
\mbox{if} \; \; \kap(j) = \infty. 
\]
\eth

\thref{main1} is equivalent to the the following theorem which will be proved 
in \S \ref{3.3}.

\bth{main1b} In the setting of \thref{main1} the quantum minors 
\eqref{minors} are expressed in terms of the generators 
$\ol{F}_{\ib, 1}, \ldots, \ol{F}_{\ib, \ell(w)}$ of the Cauchon
quantum affine space algebra associated to 
the presentation \eqref{Uwiter} 
of $\UU^-[w]$ as a CGL extension corresponding to $\ib$ by
\begin{equation}
\label{main1b-eq}
\De_{\ib, j} = (q_{\al_j}^{-1} - q_{\al_j})^{O(j)} 
\ol{F}_{\ib,\kap^{O(j)}(j)} \ldots \ol{F}_{\ib,j}, 
\; \; j \in [1, \ell(w)].
\end{equation}
\eth

The special case of this theorem for the algebras of quantum matrices 
$R_q[M_{m,n}]$ is due to Cauchon \cite{Ca2}. Given 
$m,n \in \Zset_+$, let $\g := \sl_{m+n}$ and 
$w:=w_{m,n} \in S_{m+n}$ for 
$w_{m,n} = c^m$ and $c := (1 \, 2 \, \ldots \, m+n)$. The
algebra $R_q[M_{m,n}]$ is isomorphic to the algebras 
$\UU^\pm[w_{m,n}]$ by
\cite[Proposition 2.1.1]{MC} and \cite[Lemma 4.1]{Y6}. 
In this case by \cite[Lemma 4.3]{Y6} the elements 
$b^{\vpi_\al}_{y,w_{m,n}} \in \UU^\pm[w_{m,n}]$ correspond 
(under this isomorphism) to scalar multiples of quantum minors of $R_q[M_{m,n}]$ 
for all $\al \in \Pi$, $y \in S_{m+n}^{\leq w_{m,n}}$. 
In particular, the elements 
$\De_{\ib, 1}, \ldots, \De_{\ib, mn} \in \UU^\pm[w_{m,n}]$
correspond to scalar multiples of quantum minors of $R_q[M_{m,n}]$
for all reduced words $\ib$ of $w_{m,n}$.
\subsection{Leading terms of quantum minors}
\label{3.2}
There are several different ways to construct iterated 
Ore extensions associated to the algebras $\UU^-[w]$, 
by adjoining root vectors in different order. Passing 
from one to the other will play a major role in our proof 
of \thref{main1b} in \S \ref{3.3}. 
In \S \ref{3.2}--\ref{3.3} we examine these iterated 
Ore extensions and prove a leading term result for 
the elements $\De_{\ib, j}$. 

For a reduced word $\ib=(\al_1, \ldots, \al_l)$ for $w \in W$ and
$j, k \in [1,l]$ denote by 
\[
\UU^-[w]_{\ib, [j,k]} \; \; 
\mbox{the subalgebra of $\UU^-[w]$ generated by}
\; \; F_{\be_m} \; \; \mbox{for} \; \; 
j \leq m \leq k
\] 
in terms of the notation from eq. \eqref{rootv}.
One easily shows that 
\[
\UU^-[w]_{\ib, [j,k]} = T_{w(\ib)_{\leq (j-1)}} 
(\UU^-[ (w(\ib)_{\leq (j-1)})^{-1} w(\ib)_{\leq k}]),
\]
for $j \leq k$, but we will not need this.  

\bpr{aux1} For all base fields $\KK$, 
$q \in \KK^*$ not a root of unity, simple 
Lie algebras $\g$, $w \in W$ of length $l$, 
reduced words $\ib$ for $w$, and $j \in [1,l]$,
we have
\begin{equation}
\label{id1}
\De_{\ib, j} = 
(q_{\al_j}^{-1} - q_{\al_j}) \De_{\ib, \kap(j)} F_{\be_j}
\mod \UU^-[w]_{\ib, [j+1, l]},
\; \; \mbox{if} \; \; \kap(j) \neq \infty
\end{equation}
and 
\begin{equation}
\label{id2}
\De_{\ib, j} = 
(q_{\al_j}^{-1} - q_{\al_j}) F_{\be_j}
\mod \UU^-[w]_{\ib, [j+1, l]}, 
\; \; \mbox{if} \; \; \kap(j)= \infty.
\end{equation}
\epr
\begin{proof}
We fix a reduced expression $\ib= (\al_1, \ldots, \al_l)$
of $w$ and denote $w_{\leq k } := w(\ib)_{ \leq k}$, $k \in [0,l]$, 
cf. \eqref{leq1}.
Recall that $\tau$, given by \eqref{tau}, 
is an algebra antiautomorphism of $\UU_q(\g)$ 
and $T_w$ is an algebra automorphism of $\UU_q(\g)$ for all $w \in W$.
The algebra $\tau T_{w_{\leq (k-1)}}(\UU^{\al_k})$ is (anti)isomorphic 
to $\UU_{q_{\al_k}}(\sl_2)$ for all $k \in [1, l]$, see \eqref{Ual}.
 
Let $1 \leq k < j$. 
Consider the $\tau T_{w_{\leq (k-1)}}(\UU^{\al_k})$-submodule
of $V(\vpi_{\al_j})$
generated by 
\[
V(\vpi_{\al_j})_{ w_{ \leq (j-1) } \vpi_{\al_j} } 
= \KK T_{w_{\leq (j-1)}} v_{\vpi_{\al_j}}
= \KK T_{w^{-1}_{\leq (j-1)}}^{-1} v_{\vpi_{\al_j}}.
\]
It is irreducible since 
\begin{align}
\label{kj1}
& (\tau F_{\be_k}) \left( T_{w^{-1}_{\leq (j-1)}}^{-1} v_{\vpi_{\al_j}} \right)= 
\left( \tau(T_{w_{\leq (k-1)}} F_{\al_k}) \right) 
\left( T_{w^{-1}_{\leq (j-1)}}^{-1} v_{\vpi_{\al_j}} \right) 
\\
= & \left( T_{w^{-1}_{\leq (k-1)}}^{-1} F_{\al_k} \right) 
\left( T_{w^{-1}_{\leq (j-1)}}^{-1} v_{\vpi_{\al_j}} \right) = 
T_{w^{-1}_{\leq (j-1)}}^{-1}
\left( (T_{\al_{j-1}} \ldots T_{\al_k} (F_{\al_k}))
v_{\vpi_{\al_j}}
\right) = 0,
\nn
\end{align}
cf. \eqref{compatible} and \eqref{tau-ident}. In the last equation 
we used that $-{s_{\al_{j-1}} \ldots s_{\al_k} (\al_k)} \in \QQ^+$ and
$T_{\al_{j-1}} \ldots T_{\al_k} (F_{\al_k}) \in \UU_q(\g)_{s_{\al_{j-1}} \ldots s_{\al_k} (\al_k)}$.
Therefore there exists a splitting of $\tau T_{w_{\leq (k-1)}}(\UU^{\al_k})$-modules
\[
V(\vpi_{\al_j}) = \left( \tau T_{w_{\leq (k-1)}}(\UU^{\al_k}) \right)
V(\vpi_{\al_j})_{ w_{ \leq (j-1) } \vpi_{\al_j} } \oplus V_k
\]
such that $V_k$ is also $\UU^0$-stable. From this and eq. \eqref{kj1}
it follows that 
\begin{equation}
\label{i1}
\lcor \xi_{w_{ \leq (j-1) }, \vpi_{\al_j} } ,  (\tau E_{\be_k}) v \rcor = 0, \; \; 
\forall \, 
v \in V(\vpi_{\al_j}), 1 \leq k < j,
\end{equation}  
recall \eqref{xi-w}.

Next, we consider the $\tau T_{w_{\leq (j-1)}}(\UU^{\al_j})$-submodule
of $V(\vpi_{\al_j})$ generated by $T_{w^{-1}_{\leq (j-1)}}^{-1} v_{\vpi_{\al_j}}$.
Using \eqref{sl2-braid}--\eqref{sl2-EF}, we obtain:
\begin{align}
\label{kj2}
&(\tau E_{\be_j}) \left( T_{w^{-1}_{\leq j}}^{-1} v_{\vpi_{\al_j}} \right)
= \left( \tau( T_{w_{\leq (j-1)}} E_{\al_j}) \right) 
\left( T_{w^{-1}_{\leq j}}^{-1} v_{\vpi_{\al_j}} \right)
\\
\nn
= &
\left( T_{w^{-1}_{\leq (j-1)}}^{-1} E_{\al_j} \right) 
\left( T_{w^{-1}_{\leq j}}^{-1} v_{\vpi_{\al_j}} \right)
=
T_{w^{-1}_{\leq (j-1)}}^{-1}
\left( E_{\al_j} T_{\al_j}^{-1} v_{\vpi_{\al_j} }
\right) = 
T_{w^{-1}_{\leq (j-1)}}^{-1} v_{\vpi_{\al_j}}.
\end{align}
Analogously one shows that
\begin{multline*}
( \tau E_{\be_j} ) 
\left( T_{w^{-1}_{\leq (j-1)}}^{-1} v_{\vpi_{\al_j}} \right) = 0
\; \; 
\mbox{and} 
\\
\left( \tau( T_{w_{\leq (j-1)}} K_{\al_j})  \right) 
\left( T_{w^{-1}_{\leq (j-1)}}^{-1} v_{\vpi_{\al_j}} \right) =
q_{\al_j}^{-1} \left( T_{w^{-1}_{\leq (j-1)}}^{-1} v_{\vpi_{\al_j}} \right).
\end{multline*}
Therefore 
\[
\left( \tau T_{w_{\leq (j-1)}}(\UU^{\al_j}) \right) 
T_{w^{-1}_{\leq (j-1)}}^{-1} v_{\vpi_{\al_j}} = 
\KK T_{w^{-1}_{\leq (j-1)}}^{-1} v_{\vpi_{\al_j}}
\oplus 
\KK T_{w^{-1}_{\leq j}}^{-1} v_{\vpi_{\al_j}}.
\]
Using this, the complete reducibility of finite dimensional 
type one $\UU^\al$-modules, and eq. \eqref{kj2}, we obtain:
\begin{align}
\label{i2}
&\lcor \xi_{w_{ \leq (j-1) }, \vpi_{\al_j} } ,  (\tau E_{\be_j}) v \rcor = 
\lcor \xi_{w_{ \leq j }, \vpi_{\al_j} }, v \rcor \; \; \mbox{and}
\\
&\lcor \xi_{w_{ \leq (j-1) }, \vpi_{\al_j} } ,  ( \tau E_{\be_j})^m v \rcor = 0, 
\; \; \forall \, 
v \in V(\vpi_{\al_j}), m >1.
\nn
\end{align}

In a similar way one proves that for all $j < k \leq \min \{l, \kap(j)- 1\}$
\[
\left( \tau T_{w_{\leq (k-1)}}(\UU^{\al_k}) \right) 
T_{w^{-1}_{\leq (k-1)}}^{-1} v_{\vpi_{\al_j}} = 
\KK T_{w^{-1}_{\leq (k-1)}}^{-1} v_{\vpi_{\al_j}}  \; \; 
\mbox{and} \; \;
T_{w^{-1}_{\leq (k-1)}}^{-1} v_{\vpi_{\al_j}} =
T_{w^{-1}_{\leq k}}^{-1} v_{\vpi_{\al_j}}.
\]
From this one obtains that 
\begin{multline}
\label{i3}
\lcor \xi_{w_{ \leq (k-1) }, \vpi_{\al_j} } ,  (\tau E_{\be_k}) v \rcor = 0
\; \; 
\mbox{and}
\\
\xi_{w_{ \leq (k-1) }, \vpi_{\al_j} } =
\xi_{w_{ \leq  k}, \vpi_{\al_j} },
\; \; \forall \, 
v \in V(\vpi_{\al_j}), 
j < k \leq \min \{l, \kap(j)- 1\}.
\end{multline}

Eq. \eqref{id1} is deduced from from eqs. \eqref{i1}, \eqref{i2}, 
and \eqref{i3} as follows. Denote for brevity
\[
p_{k, m} := \frac{ (q_{\al_k}^{-1} - q_{\al_k})^m}
{q_{\al_k}^{m (m-1)/2} [m]_{\al_k}! }, \; \; 
k \in [1,l], m \in \Nset.
\]
Using \eqref{phi1}, \eqref{i1}, and \eqref{i2}, we obtain:
\begin{align*} 
&\De_{\ib, j} = \sum_{m_j, \ldots, m_l \in \Nset}
\Big( \prod_{k=j}^l 
p_{k, m_k} \Big) 
\lcor \xi_{w_{\leq (j-1)}, \vpi_{\al_j}}, 
(\tau E_{\be_j})^{m_j} \ldots (\tau E_{\be_l})^{m_l} T^{-1}_{w^{-1}} v_\la \rcor 
F_{\be_l}^{m_l} \ldots F_{\be_j}^{m_j}
\\
= &(q_{\al_j}^{-1} - q_{\al_j})
\sum_{m_{j+1}, \ldots, m_l \in \Nset}
\Big( \prod_{k=j+1}^l 
p_{k, m_k} \Big) 
\lcor \xi_{w_{\leq j}, \vpi_{\al_j}}, 
(\tau E_{\be_{j+1}})^{m_{j+1}} \ldots (\tau E_{\be_l})^{m_l} T^{-1}_{w^{-1}} v_\la \rcor 
\\
& \hspace{1cm} \times 
F_{\be_l}^{m_l} \ldots F_{\be_{j+1}}^{m_{j+1}} F_{\be_j}
\mod \UU^-[w]_{\ib, [j+1, l]}.
\end{align*}
If $\kap(j) \leq l$, it follows from \eqref{i3} that the right hand side 
of the last congruence is equal to
\begin{align*}
&(q_{\al_j}^{-1} - q_{\al_j})
\sum_{m_{\kap(j)}, \ldots, m_l \in \Nset}
\Big( \prod_{k=\kap(j)}^l 
p_{k, m_k} \Big) 
\lcor \xi_{w_{\leq (\kap(j)-1)}, \vpi_{\al_j}}, 
(\tau E_{\be_{\kap(j)}})^{m_{\kap(j)}} \ldots (\tau E_{\be_l})^{m_l} T^{-1}_{w^{-1}} v_\la \rcor 
\\
&\hspace{1cm} \times 
F_{\be_l}^{m_l} \ldots F_{\be_{\kap(j)}}^{m_{\kap(j)}} F_{\be_j}
=(q_{\al_j}^{-1} - q_{\al_j}) \De_{\ib, \kap(j)} F_{\be_j}.
\end{align*}
This proves eq. \eqref{id1}. 
The proof of eq. \eqref{id2} is analogous, requiring 
only a small modification of the last argument. It
is left to the reader. 
\end{proof}

Starting from a reduced word $\ib=(\al_1, \ldots, \al_l)$ 
for $w \in W$, one can
construct a presentation of $\UU^-[w]$ as an iterated 
Ore extension by adjoining the elements $F_{\be_1},$ $\ldots$, 
$F_{\be_l}$ (recall \eqref{rootv}) in the opposite order.
For all $j \in [1,l]$ we have the Ore extension presentation 
\begin{equation}
\label{Urev}
\UU^-[w]_{\ib, [j,l]} = 
\UU^-[w]_{\ib, [j+1, l]}[F_{\be_j}; \sig'_j, \delta'_j],
\end{equation}
where $\sig'_j$ and $\delta'_j$ are defined as follows.
Let $t'_j$ be an element of $\Tset^r$ such that 
\[
(t'_j)^{\be_k} = 
q^{ - \lcor \be_k, \be_j \rcor }, \; \; 
\forall \, k \geq j
\]
(cf. \eqref{Tchar} and \eqref{tj}) and $\sig'_j:=(t'_j \cdot)$ in terms of the 
restriction of the $\Tset^r$-action \eqref{torus-act} to $\UU^-[w]_{\ib, [j+1, l]}$.
The skew derivation $\delta'_j$ of $\UU^-[w]_{\ib, [j+1, l]}$
is defined by
\[
\delta'_j(x) := F_{\be_j} x - q^{ - \lcor \be_j, \ga \rcor} x F_{\be_j}, \; \; 
x \in (\UU^-[w]_{\ib, [j+1,l]})_\ga, \ga \in \QQ,
\]
cf. \eqref{deltaj}. (It follows from the Levendorskii--Soibelman 
straightening law \eqref{LS} 
that $\delta'_j$ preserves $\UU^-[w]_{\ib, [j+1,l]}$, $\sig'_l=\id$,
and $\delta'_l = 0$.) Eqs. \eqref{PBW} and \eqref{LS} imply \eqref{Urev}.
Iterating \eqref{Urev} and taking into account 
$\UU^-[w]_{\ib, [l+1, l]} = \KK$
leads to the iterated Ore extension presentation
\[
\UU^-[w] = \KK[F_{\be_l}][F_{\be_{l-1}}; \sig'_{l-1}, \de'_{l-1}]
\ldots [F_{\be_1}; \sig'_1, \de'_1],
\]
which is reverse to the presentation \eqref{Uwiter}. 
It is straightforward to show that this presentation of $\UU^-[w]$ 
is a torsion free CGL extension for the action \eqref{torus-act}.

In this framework, \prref{aux1} proves that 
$\De_{\ib, j} \in \UU^-[w]_{\ib, [j, l]}$ 
and computes its leading term as a left polynomial with respect
to the Ore extension \eqref{Urev}, for all $j \in [1,l]$, 
cf. \S \ref{2.4}.
\subsection{Proof of \thref{main1b}}
\label{3.3}
We keep the notation for $\ib$, $w$, and $l$ from the 
previous two subsections.
For $j \in [1,l]$ consider the chain of extensions
\[
\KK \subset \UU^-[w]_{\ib, [j,j]} \subset 
\UU^-[w]_{\ib, [j, j+1]} \subset \ldots \subset
\UU^-[w]_{\ib, [j, l]}.
\]
It follows from the Levendorskii--Soibelman 
straightening law \eqref{LS} and the definition of the 
$\Tset^r$-action \eqref{torus-act} that the maps $\delta_k$ and $\sig_k$ 
from \leref{Uw} (ii) preserve the subalgebra $\UU^-[w]_{\ib, [j,k-1]}$ 
of $\UU^-[w(\ib)_{\leq (k-1)}] = \UU^-[w]_{\ib, [1,k-1]}$
for all 
$1 \leq j \leq k \leq l$. Denote the restrictions 
\[
\delta_{j,k} = \delta_k |_{\UU^-[w]_{\ib, [j,k-1]}} \; \; 
\mbox{and} \; \; 
\sig_{j,k} = \sig_k |_{\UU^-[w]_{\ib, [j,k-1]}}, \; \; 
\mbox{for} \; \; 1 \leq j \leq k \leq l.
\]
\leref{Uw} (ii) implies that we have the Ore extension presentation 
\[
\UU^-[w]_{\ib, [j, k]} = 
\UU^-[w]_{\ib, [j, k-1]}[F_{\be_k}; \sig_{j,k}, \delta_{j,k}], \; \; 
\mbox{for} \; \; 1 \leq j \leq k \leq l.
\] 
Iterating those and using that $\UU^-[w]_{\ib, [j, j-1]}= \KK$,
$\sig_{j,j} = \id$, and
$\delta_{j,j} = 0$ leads to the iterated Ore extension presentation 
of $\UU^-[w]_{\ib, [j,k]}$:
\begin{equation}
\label{Ujl}
\UU^-[w]_{\ib, [j,l]} = \KK [F_{\be_j}]
[F_{\be_{j+1}}; \sig_{j, j+1}, \delta_{j, j+1}] \ldots 
[F_{\be_l}; \sig_{j, l}, \delta_{j,l}]. 
\end{equation}
It follows now from \leref{Uw} that $\UU^-[w]_{\ib, [j,k]}$
is a CGL extension with respect to the $\Tset^r$-action 
\eqref{torus-act}. Since $\{ 0 \}$ is a $\Tset^r$-prime ideal of 
$\UU^-[w]_{\ib, [j,k]}$, we can 
apply a theorem of Goodearl \cite[Theorem II.6.4]{BG}, to obtain 
that it is a strongly rational ideal, i.e.,
\begin{equation}
\label{str-r}
Z(\Fract(\UU^-[w]_{\ib, [j,l]}))^{\Tset^r} = \KK.
\end{equation}
Recall that $Z(A)$ stands for the center of an algebra $A$.
As in \S \ref{2.4}, $\Fract(A)$ denotes the division ring 
of fractions of a domain $A$. Furthermore, $(.)^{\Tset^r}$ 
refers to the fixed point subalgebra with respect 
to the action \eqref{torus-act}.

Denote by $\TT_{\ib}$ the quantum torus algebra generated 
by $\ol{F}_{\ib, 1}^{\pm 1}, \ldots, \ol{F}_{\ib, l}^{\pm}$. 
Eqs. \eqref{LS} and \eqref{final-isom} imply that 
\begin{equation}
\label{Fcomm}
\ol{F}_{\ib, j} \ol{F}_{\ib, k} = 
q^{ \lcor \be_j, \be_k \rcor} \ol{F}_{\ib, k} \ol{F}_{\ib, j}, \; \; 
\forall \, 1 \leq j < k \leq l.
\end{equation}
For $j, k \in [1,l]$
denote by $\TT_{\ib, [j,k]}$ the quantum subtorus 
of $\TT_\ib$ generated by $\ol{F}_{\ib, m}^{\pm 1}$ for 
$j \leq m \leq k$.

Using that
\[
\delta_k (F_{\be_j}) \in \UU^-[w]_{\ib, [k+1, j-1]},
\]
by a simple induction argument one proves the following 
lemma:
\ble{aux2} In the above setting, the following hold  for all $j \in [1,l]$:  

(i) $F_{\be_j} - \ol{F}_{\ib, j} \in \TT_{\ib, [j+1,l]}$. 

(ii) The generators for the Cauchon quantum affine space algebra 
associated to the iterated Ore extension presentation \eqref{Ujl}
of $\UU^-[w]_{\ib, [j,l]}$ are precisely the elements 
$\ol{F}_{\ib, j}$, $\ldots$, $\ol{F}_{\ib, l}$, 
recall \S \ref{2.4}.
\ele

The lemma implies that 
$\UU^-[w]_{\ib, [j,l]} \subset
\TT_{\ib, [j,l]} \subset \Fract(\UU^-[w]_{\ib, [j,l]})$. 
Therefore the strong rationality result \eqref{str-r} gives that 
\begin{equation}
\label{str-r0}
Z(\TT_{\ib, [j,l]})_0 = \KK,
\end{equation}  
where $(.)_0$ refers to the $0$-component with 
respect to the $\QQ$-grading induced from the 
grading of $\UU_q(\g)$.

Next we apply a theorem of Berenstein and Zelevinsky 
\cite[Theorem 10.1]{BZ}, to obtain that 
there exist integers $n_{j k } \in \Zset$ ($1 \leq j < k \leq n$) such that 
\[
e^{\vpi_{\al_j}}_{w(\ib)_{\leq (j-1)}} e^{\vpi_{\al_k}}_{w(\ib)_{\leq (k-1)}}
= q^{n_{jk}}
e^{\vpi_{\al_k}}_{w(\ib)_{\leq (k-1)}} e^{\vpi_{\al_j}}_{w(\ib)_{\leq (j-1)}},
\; \; 
\forall \, 
1 \leq j < k \leq l.  
\]  
(The setting of \cite{BZ} is for $\KK= \Qset(q)$, but the proof
of Theorem 10.1 in \cite{BZ} only uses the $R$-matrix commutation 
relations in $R_q[G]$ and the left and right actions of $\UU_q(\g)$ on 
$R_q[G]$, which work for all fields $\KK$ and $q \in \KK^*$ 
not a root of unity.) Moreover, the $R$-matrix commutation relations 
in $R_q[G]$ (see e.g. \cite[Theorem I.8.15]{BG}) imply that 
\[
e^\la_w c^{\la'}_{\xi'} = q^{-\lcor \la,\la' + w^{-1}\mu'\rcor} 
c^{\la'}_{\xi'} e^\la_w \mod Q(w)^+, \; \; 
\forall \la, \la' \in \PP^+,
\mu \in \PP, \xi' \in V(\la')_{\mu'}.
\]
Using \eqref{minors} and the fact that the maps
$\phi_w \colon R_w^0 \to \UU^-[w]$ are antihomomorphisms 
by \thref{Class} (i), 
we obtain 
\begin{equation}
\label{d-com}
\De_{\ib, j} \De_{\ib, k} = q^{n'_{jk}} \De_{\ib, k} \De_{\ib, j}, 
\; \; 
\forall \, 1 \leq j < k \leq l
\end{equation}
for some $n'_{jk} \in \Zset$.
\medskip
\\
\noindent
{\em{Proof of \thref{main1b}}}. By \leref{aux2} (ii)
\[
\UU^-[w]_{\ib, [j,l]} \subseteq \TT_{\ib, [j,l]}, \; \; 
\forall \, j \in [1,l].
\]
Combining this, \prref{aux1}, and \leref{aux2} (i), 
we obtain
\begin{equation}
\label{id3}
\De_{\ib, j} = 
(q_{\al_j}^{-1} - q_{\al_j}) \De_{\ib, \kap(j)} \ol{F}_{\ib, j}
\mod \TT_{\ib, [j+1, l]},
\; \; \mbox{if} \; \; \kap(j) \leq l
\end{equation}
and 
\begin{equation}
\label{id4}
\De_{\ib, j} = 
(q_{\al_j}^{-1} - q_{\al_j}) \ol{F}_{\ib, j}
\mod \TT_{\ib, [j+1, l]}, 
\; \; \mbox{if} \; \; \kap(j)= \infty.
\end{equation}

We prove eq. \eqref{main1b-eq} by induction on $j$, from $l$ to $1$. 
By \eqref{id4}, $\De_{\ib, l} - (q_{\al_l}^{-1} - q_{\al_l}) \ol{F}_{\ib, l} \in \KK$.
Since $\De_{\ib, l}$ is a homogeneous element of nonzero degree (equal to $\be_l$),
this implies \eqref{main1b-eq} for $j=l$. 

Now assume that for some $j \in [1, l-1]$ 
\begin{equation}
\label{ind-assumpt}
\De_{\ib, k} = (q_{\al_k}^{-1} - q_{\al_k})^{O(k)} 
\ol{F}_{\ib,\kap^{O(k)}(k)} \ldots \ol{F}_{\ib,k}
\; \; \mbox{for all} \; \; k \in [j+1, l].
\end{equation}
If 
\begin{equation}
\label{ind-st}
\De_{\ib, j} = (q_{\al_j}^{-1} - q_{\al_j})^{O(j)} 
\ol{F}_{\ib,\kap^{O(j)}(j)} \ldots \ol{F}_{\ib,j}, 
\end{equation}
then we are done with the inductive step. 
Assume the opposite, that \eqref{ind-st} is 
not satisfied. Combining the inductive hypothesis with
\eqref{id3} and 
\eqref{id4} (whichever applies for the particular $j$), 
we get that 
\begin{equation}
\label{l5}
\De_{\ib, j} -
(q_{\al_j}^{-1} - q_{\al_j})^{O(j)} 
\ol{F}_{\ib, \kap^{O(j)}(j)} \ldots \ol{F}_{\ib, j}
\in \TT_{\ib, [j+1, l]}.
\end{equation}
It follows from eqs. \eqref{Fcomm},
\eqref{d-com}, and \eqref{ind-assumpt}, 
that 
\[
\De_{\ib,j} \ol{F}_{\ib, k} = q^{m_k} \ol{F}_{\ib, k} 
\De_{\ib,j}, \; \; 
\forall \, k =j+1, \ldots, l  
\] 
for some $m_{j+1}, \ldots, m_l \in \Zset$. Quantum tori
have bases consisting of Laurent monomials in their 
generators. By comparing the coefficients of 
$\ol{F}_{\ib, \kap^{O(j)}(j)} \ldots \ol{F}_{\ib, j} \ol{F}_{\ib, k}$
in the two sides of the above equality 
and using \eqref{l5}, we get that 
\[
(\ol{F}_{\ib, \kap^{O(j)}(j)} \ldots \ol{F}_{\ib, j})
\ol{F}_{\ib, k} = q^{m_k} \ol{F}_{\ib, k}
(\ol{F}_{\ib, \kap^{O(j)}(j)} \ldots \ol{F}_{\ib, j}), \; \; 
\forall \, k =j+1, \ldots, l  
\]
for the same collection of integers $m_{j+1}, \ldots, m_l$.
From the last two equalities it follows that 
\[
y : = (\ol{F}_{\ib,\kap^{O(j)}(j)} \ldots \ol{F}_{\ib,j})^{-1} \De_{\ib, j}
\] 
commutes with $\ol{F}_{\ib, j+1}, \ldots, \ol{F}_{\ib, l}$:
\begin{equation}
\label{commute5}
y \ol{F}_{\ib, k} = \ol{F}_{\ib, k} y, \; \; 
\forall \, k =j+1, \ldots, l.  
\end{equation}
Since 
\eqref{ind-st} is not satisfied, \eqref{l5} implies that 
\begin{equation}
\label{last}
y = (q_{\al_j}^{-1} - q_{\al_j}) + y' \ol{F}_{\ib, j}^{\, -1}
\; \; 
\mbox{for some} \; \; y' \in \TT_{\ib, [j+1, l]} \backslash \{ 0\}. 
\end{equation}
But $y$ commutes with itself and by \eqref{commute5} it commutes
with $y' \neq 0$. Thus 
$y$ also commutes with $\ol{F}_{\ib, j}$. Combining this with 
\eqref{commute5} leads to the fact that $y$ belongs to the center 
of $\TT_{\ib, [j,l]}$. Since 
$\De_{\ib, j}$ is a homogeneous element of $\UU^-[w]$ with respect 
to its $\QQ$-grading, \eqref{l5} implies
\[
y \in Z(\TT_{\ib, [j,l]})_0.
\] 
At the same time $y \notin \KK$ by \eqref{last}, which contradicts with 
the strong rationality result \eqref{str-r0}. Thus
\eqref{ind-st} holds. This completes the proofs 
of the inductive step and the theorem.
\qed  
\sectionnew{Unification of the two approaches to 
$\TSpec \UU^-[w]$} 
\label{two-appr}
\subsection{Solutions of two questions of Cauchon and M\'eriaux}
\label{4.1}
In this section we establish a relationship between the 
representation theoretic and ring theoretic
approaches to the prime spectra of the quantum Schubert 
cell algebras $\UU^-[w]$, see \S \ref{2.3} and \S \ref{2.4}. 
\thref{main2-ind} explicitly describes the behavior of all 
$\Tset^r$-prime ideals $I_w(y)$ of the algebras $\UU^-[w]$ 
from \thref{Class} under the iterations of Cauchon's deleting derivation 
construction, recall \prref{ind}. In \thref{main2} we describe explicitly 
the Cauchon diagrams of all ideals $I_w(y)$ and use this to resolve 
\cite[Question 5.3.3]{MC} of Cauchon and Meriaux. We use the combination 
of Theorems \ref{tClass} and \ref{tmain2} to give a new,
independent proof of the classification result in \thref{CM} 
of Cauchon and Meriaux. Finally, we also 
settle \cite[Question 5.3.2]{MC} of Cauchon and Meriaux, 
solving the containment problem for the ideals in the 
classification of \thref{CM}, see \thref{co2}.

\bth{main2} Assume that $\KK$ is an arbitrary base field, 
$q \in \KK^*$ is not a root of unity, $\g$ is a simple 
Lie algebras $\g$, $w$ is a Weyl group element, 
and $\ib$ is a reduced word for $w$. Then for all Weyl group 
elements $y \leq w$ the Cauchon diagram of the 
$\Tset^r$-prime ideal $I_w(y)$ (see \thref{Class} (ii))
for the presentation \eqref{Uwiter} of $\UU^-[w]$
is precisely the index set of the left positive 
subword of $\ib$ whose total product is $y$ 
\[
\CD(I_w(y)) = \LP_\ib(y),
\] 
recall \S 2.2 and 2.4 for definitions.
\eth 
\bre{new} \thref{main2} gives a new, independent proof of 
\thref{CM} of Cauchon and M\'eriaux \cite{MC}. By \thref{Class} (ii)
\[
\TSpec \UU^-[w] = \{ I_w(y) \mid y \in W^{\leq w} \}.
\]
Since $\CD(I_w(y)) = \LP_\ib(y)$ by \thref{main2} we have
\[
\TSpec \UU^-[w] = \{ J_{\LP_\ib(y)} \mid y \in W^{\leq w} \},
\] 
which is the statement of \thref{CM}, recall \eqref{JD}.
\ere

The following theorem is an immediate consequence of \thref{main2}.
It settles Question 5.3.3 of Cauchon and M\'eriaux \cite{MC}.

\bth{co1} For all base fields $\KK$, $q \in \KK^*$ not a root 
of unity, simple Lie algebras $\g$, Weyl group elements $w$, 
and reduced words $\ib$ for $w$,
\begin{equation}
\label{eq-id}
I_w(y) = J_{\LP_\ib(y)}, \; \; \forall \, y \in W^{\leq w}
\end{equation}
(recall \eqref{JD}),
i.e., the the classifications of $\TSpec \UU^-[w]$ of 
Cauchon--M\'eriaux \cite{MC} from \thref{CM} 
and Yakimov \cite{Y1} from \thref{Class} coincide.
\eth

Finally, the next theorem answers Question 5.3.2 of 
Cauchon and M\'eriaux \cite{MC}.
\bth{co2} For all base fields $\KK$, $q \in \KK^*$ not a root 
of unity, simple Lie algebras $\g$, Weyl group elements $w$, 
and reduced words $\ib$ for $w$, the map
\[
y \in W^{\leq w} \mt J_{\LP_\ib(y)} \in \TSpec \UU^-[w], 
\quad y \in W^{\leq w},
\]
is an isomorphism of posets with respect to the Bruhat order and 
inclusion of ideals.
\eth
\begin{proof} \thref{co2} follows from \thref{Class} (iii) and eq. \eqref{eq-id}.
\end{proof}

Our proof of \thref{main2} is based on a result, which gives
a full picture of the 
behavior of the ideals $I_w(y)$ from \thref{Class} (i) 
under the deleting derivation procedure from \S \ref{2.4}.
Recall the definition \eqref{hI} of leading part $\lt(J)$ of an
ideal of an Ore extension.
According to \prref{ind}, Cauchon's method relies on 
taking leading parts or contractions of ideals in 
CGL extensions. Assume that $\ib=(\al_1, \ldots, \al_l)$ 
is a reduced word for $w \in W$. Then  
\begin{equation}
\label{Ore-1}
w(\ib)_{\leq (l-1)} = w s_{\al_l}.
\end{equation}
\leref{Uw} (i)--(ii) implies that 
\begin{equation} 
\label{Ore-2} 
\UU^-[w s_{\al_l}] = \UU^-[w(\ib)_{\leq (l-1)}] 
\subset \UU^-[w] \; \; 
\mbox{and} \; \;  
\UU^-[w]
= \UU^-[w s_{\al_l}][F_{\be_l}; \sig_l, \delta_l],
\end{equation}
where $\sig_l$ and $\delta_l$ are the automorphism 
and left $\sig_l$-skew derivation of 
$\UU^-[w(\ib)_{\leq (l-1)}]$ from \leref{Uw} (ii).
 We have:

\bth{main2-ind} Assume that $\KK$ is an arbitrary base field, 
$q \in \KK^*$ is not a root of unity, $\g$ is a simple 
Lie algebra, $w \in W$ is a Weyl group element of length $l$, 
and $\ib=(\al_1, \ldots, \al_l)$ is 
a reduced word for $w$. Then the following hold for all $y \in W^{\leq w}$:

(i) If $l \notin \LP_\ib(y)$, then $\lt(I_w(y)) = I_{w s_{\al_l}}(y)$, 
where the leading part of $I_w(y)$ (cf. \eqref{hI}) is computed with respect to the 
Ore extension $\UU^-[w] = \UU^-[w s_{\al_l}][F_{\be_l}; \sig_l, \delta_l]$, 
cf. \eqref{Ore-2}.

(ii) If $l \in \LP_\ib(y)$, then $I_w(y) \cap \UU^-[ w s_{\al_l}] 
= I_{w s_{\al_l}}( y s_{\al_l})$.   
\eth
We prove \thref{main2} using \thref{main2-ind} in this subsection.
We establish \thref{main2-ind} in \S \ref{4.2}--\ref{4.3}.
Before we proceed with the proof of \thref{main2}, we 
prove an auxiliary lemma.
\ble{lem} If, in the setting of \thref{main2-ind}, $y \in W^{\leq w}$ 
is such that $l \in \LP_\ib(y)$, then 
\begin{equation}
\label{incl5}
T_{w s_{\al_l}} v_{\vpi_{\al_l}} \notin \UU^- T_y v_{\vpi_{\al_l}}.
\end{equation}
\ele
\begin{proof} The similar statement that 
$T_{w s_{\al_l}} v_\la \notin \UU^- T_y v_\la$ for 
$\la \in \sum_{\al \in \Pi} \Zset_+ \vpi_\al$ follows from
\cite[Lemma 4.4.5]{J} and the fact that $y \not\leq w s_{\al_l}$,
which is easy to show. The last lemma is not applicable in
our case, but we use some ideas of its proof.

We argue by induction on $l= \ell(w)$. If $l=1$, then 
$T_{w s_{\al_l}} v_{\vpi_{\al_l}} = v_{\vpi_{\al_l}}$ and the statement 
is true since $y (\vpi_{\al_l}) < \vpi_{\al_l}$.
Assume the validity of the lemma for length $l-1$. 

Let $y \leq w \in W$ and $\ib$ be as in the statement of the lemma. Assume 
that \eqref{incl5} does not hold, i.e.,
\begin{equation}
\label{incl6}
T_{w s_{\al_l}} v_{\vpi_{\al_l}} \in \UU^- T_y v_{\vpi_{\al_l}}.
\end{equation}
We consider two cases: (A) $1 \in \LP_\ib(y)$ and (B) $1 \notin \LP_\ib(y)$.
Note that 
\[
\ib'' := (\al_2, \ldots, \al_l) \; \; 
\mbox{is a reduced word for} \; \; s_{\al_1} w.
\]

Case (A) $1 \in \LP_\ib(y)$. Using the left positivity 
of the index set $\LP_\ib(y)$, we obtain
\begin{equation}
\label{eqq4}
y= s_{\al_1} w(\ib) ^{\LP_\ib(y)}_{> 1}
> w(\ib)^{\LP_\ib(y)}_{> 1} = s_{\al_1} y.
\end{equation}
Moreover, we have $s_{\al_1} y \leq s_{\al_1} w$ and 
$\LP_{\ib''}( s_{\al_1} y) = \LP_{\ib} (y) \backslash \{ 1 \}$.  
Recall the definition \eqref{Ual} of the subalgebras $\UU^\al$ of $\UU_q(\g)$,
$\al \in \Pi$.
Eq. \eqref{incl6}, \eqref{eqq4} and \cite[Lemma 4.4.3 (iii)--(iv)]{J} imply
\[
T_{s_{\al_1} w s_{\al_l}} v_{\vpi_{\al_l}} 
\in \UU^{\al_1} T_{w s_{\al_l}} v_{\vpi_{\al_l}} \subseteq
\UU^{\al_1} \UU^- T_{y} v_{\vpi_{\al_l}} = 
\UU^- \UU^{\al_1} T_{y} v_{\vpi_{\al_l}} = 
\UU^- T_{s_{\al_1} y} v_{\vpi_{\al_l}},
\]
which contradicts with the induction assumption for the triple
$(s_{\al_1} y, s_{\al_1} w, \ib'')$.

Case (B) $1 \notin \LP_\ib(y)$. The argument in this case is 
similar to the previous one.
From the left positivity 
of the index set $\LP_\ib(y)$ we have
\begin{equation}
\label{eqq5}
s_{\al_1} y= s_{\al_1} w(\ib) ^{\LP_\ib(y)}_{> 1}
> w(\ib)^{\LP_\ib(y)}_{> 1} = y.
\end{equation}
Furthermore, $y < s_{\al_1} w$ and $\LP_{\ib''}(y) = \LP_{\ib}(y)$. 
Eqs. \eqref{incl6}, \eqref{eqq5} and 
\cite[Lemma 4.4.3 (iii)--(iv)]{J} imply
\[
T_{s_{\al_1} w s_{\al_l}} v_{\vpi_{\al_l}} 
\in \UU^{\al_1} T_{w s_{\al_l}} v_{\vpi_{\al_l}} \subseteq
\UU^{\al_1} \UU^- T_{y} v_{\vpi_{\al_l}} = 
\UU^- \UU^{\al_1} T_{y} v_{\vpi_{\al_l}} 
= \UU^- T_{y} v_{\vpi_{\al_l}}.
\]
This contradicts with the induction assumption for the triple
$(y, s_{\al_1} w, \ib'')$.

We reached a contradiction in both cases. Thus \eqref{incl6} is 
incorrect, which completes the proof of the lemma.
\end{proof}
\noindent
{\em{Proof of \thref{main2}}}. 
We prove \thref{main2} by induction on the length $l = \ell(w)$.
The case $\ell(w) =0$ is trivial.
Assume the validity of the statement of the theorem for length $l-1$.  

Fix $w \in W$ and a reduced word 
$\ib=(\al_1, \ldots, \al_l)$ for it.
Denote the reduced word
\[
\ib' := (\al_1, \ldots, \al_{l-1})
\]
for $w s_{\al_l}$. In the setting of \S \ref{2.4}, $\ol{x}_l = x_l$.
\thref{main1} implies that 
\[
F_{\be_l} = p_l \De_{\ib, l} = p_l b^{\vpi_{\al_l}}_{w s_{\al_l}, w}
\]
for some $p_l \in \KK^*$.
Let $y \in W^{\leq w}$. We have two cases: 
(1) $l \notin \LP_{\ib}(y)$ and (2) $l \in \LP_{\ib}(y)$.   
For brevity, in this proof we set
\[
D:= \LP_{\ib}(y).
\]

Case (1) $l \notin D$. In this case 
$w(\ib)^D_{> j} = (w s_{\al_l})(\ib')^D_{> j}$ for all $j \in [0, l-1]$. 
Taking into account \eqref{l-positive}, one sees that
$D \subseteq [1,l-1]$ is the index set of a left positive 
subword of $\ib'$. Therefore $y = (w s_{\al_l})^D < w s_{\al_l}$ 
and $\LP_{y}(\ib') = D$. The inductive assumption applied to 
$y \leq w s_{\al_l}$ implies
\begin{equation}
\label{CDind1}
\CD(I_{w s_{\al_l}} (y)) = D. 
\end{equation}
Recall from \S \ref{2.3} that 
$b^{\vpi_{\al_l}}_{w s_{\al_l}, w} \notin I_w(w s_{\al_l})$,
see \cite[Theorem 3.1 (b)]{Y6} for a proof. Thus 
$F_{\be_l} = p_l b^{\vpi_{\al_l}}_{w s_{\al_l}, w} \notin I_w(w s_{\al_l})$, 
because $p_l \in \KK^*$. \thref{Class} (ii) implies that 
$I_w(y) \subseteq I_w(w s_{\al_l})$. Therefore 
$F_{\be_l} \notin I_w(y)$. Now we are in the situation of part 
(i) of \prref{ind} with respect to the iterated Ore extension
from \eqref{Uwiter} and the ideal $J=I_w(y)$.
By \thref{main2-ind} (i), $\lt(I_w(y))= I_{w s_{\al_l}}(y)$ 
and from \prref{ind} (i) we obtain that 
$\CD (I_w(y)) = \CD(I_{w s_{\al_l}}(y))$. It follows from this 
and eq. \eqref{CDind1} that in the first case $\CD(I_w(y)) = D = \LP_\ib(y)$.

Case (2) $l \in D$. Denote $D' = D \backslash \{ l \}$. Since 
$D = \LP_\ib(y)$ we have $s_{\al_j}w(\ib)_{>j}^D > w(\ib)_{>j}^D$, 
$\forall j \in [1, l]$. 
Moreover, $w(\ib)_{>j}^D = (w s_{\al_l})(\ib')_{>j}^{D'} s_{\al_l}$
and $\ell(w(\ib)_{>j}^D) = \ell((ws_{\al_l})(\ib')_{>j}^{D'})+ 1$. This implies that 
$s_{\al_j} \big( (ws_{\al_l})(\ib')_{>j}^{D'} \big) > (ws_{\al_l})(\ib')_{>j}^{D'}$, 
$\forall j \in [1, l-1]$.
Therefore $D'$ is the index set of a left positive subword of $\ib'$.
Because $y = w(\ib)^D = (ws_{\al_l})(\ib')^{D'} s_{\al_l}$, we have
$D' = \LP_{\ib'}(y s_{\al_l})$. The inductive assumption, 
applied to $y s_{\al_l} \leq w s_{\al_l}$, implies 
\begin{equation}
\label{CDind2}
\CD(I_{w s_{\al_l}} (y s_{\al_l})) = D' = D \backslash \{ l \}. 
\end{equation}
\leref{lem} asserts that 
$T_{w s_{\al_l}} v_{\vpi_{\al_l}} \notin \UU^- T_y v_{\vpi_{\al_l}}$, so 
$\xi_{w s_{\al_l}, \vpi_{\al_l}} \in (\UU^- T_y v_{\vpi_{\al_l}})^\perp$ and 
$F_{\be_l} = p_l b^{\vpi_{\al_l}}_{w s_{\al_l}, w} \in I_w(y)$.
We are in the situation of part (ii) of \prref{ind} with respect 
to the iterated Ore extension from \eqref{Uwiter} and the ideal $J=I_w(y)$.
\thref{main2-ind} (ii) implies  
$I_w(y) \cap \UU^-[w s_{\al_l}]= I_{w s_{\al_l}}(y s_{\al_l})$. 
It follows from \prref{ind} (i) and eq. \eqref{CDind2} that 
$\CD (I_w(y)) = \CD(I_{w s_{\al_l}}(y s_{\al_l})) \sqcup \{l\}
= D' \sqcup \{ l \} = \LP_\ib(y)$.
\qed
\subsection{Proof of the first part of \thref{main2-ind}}
\label{4.2}
Recall that in the setting of \thref{main2-ind} we have the Ore extension 
$\UU^-[w] = \UU^-[w s_{\al_l}][F_{\be_l}; \sig_l, \delta_l]$
from \eqref{Ore-2}. We will prove the first part of 
\thref{main2-ind} by showing that the leading part $\lt(I_w(y))$ of the 
ideal $I_w(y)$ with respect to this Ore extension contains 
the ideal $I_{w s_{\al_l}}(y)$. We will then compare the  
Gelfand--Kirillov dimensions of the quotients 
$\UU^-[w]/I_w(y)$ and $\UU^-[w s_{\al_l}]/\lt(I_w(y))$
using results of \cite{Y5} and \prref{ind} (i) 
to show that the leading part $\lt(I_w(y))$ is precisely 
$I_{w s_{\al_l}}(y)$. The first part of this argument is based on: 

\bpr{aux2} For all base fields $\KK$, $q \in \KK^*$ not a root of unity, 
Weyl group elements $w \in W$, reduced words 
$\ib = ( \al_1, \ldots, \al_l)$ for $w$,  
$\la \in \PP^+$, and $\xi \in V(\la)^*$,
we have
\[
\phi_w( c_\xi^\la e_w^{-\la} ) 
- 
(q_{\al_l}^{-1} - q_{\al_l})^Nq_{\al_l}^{- N (N-1)/2}
F_{\be_l}^{N} 
\phi_{w s_{\al_l}}( c_\xi^\la e_{w s_{\al_l}}^{-\la} ) 
\in \sum_{m=0}^{N -1 } 
F_{\be_l}^m \UU^-[w s_{\al_l}] ,
\]
where $N:= \lcor \la, \al_l\spcheck \rcor$,
(recall \eqref{rootv}, \eqref{phi}, and \eqref{Ore-2}).
\epr

\prref{aux2} computes the leading term of $\phi_w(c_\xi^\la e^{-\la}_w)$ 
written as a right polynomial in $F_{\be_l}$ with coefficients in 
$\UU^-[w s_{\al_l}]$ (with respect to the Ore extension
\eqref{Ore-2}) if this polynomial has degree equal 
to $\lcor \la, \al_l\spcheck \rcor$, which is the highest 
expected degree. This proposition can be viewed as a dual
result to \prref{aux1}. 
\medskip
\\
\noindent
{\em{Proof of \prref{aux2}.}}
Set 
\[
w':= w s_{\al_l} = w(\ib)_{\leq (l-1)}.
\]
Recall \eqref{Ual}.
The vector $v_\la$ is a highest weight vector for $\UU^{\al_l}$ of highest 
weight $N \vpi_{\al_l}$. Eqs. \eqref{sl2-braid} and \eqref{sl2-EF}
imply
\[
E_{\al_l}^N T_\al^{-1} v_\la = \frac{1}{[N]_{\al_l}!} 
E_{\al_l}^N F_{\al_l}^N v_\la = [N]_{\al_l}! v_{\al} \; \; 
\mbox{and} \; \;  
E_{\al_l}^m T_\al^{-1} v_\la = 0, \; \; \forall \, m > N.
\]
Therefore
\[
(\tau E_{\be_l})^N T_{w^{-1}}^{-1} v_\la = 
\left( T_{(w')^{-1}}^{-1} (E_{\al_l}^N) \right)
\left( T_{(w')^{-1}}^{-1} T_\al^{-1} v_\la \right) =
T_{(w')^{-1}}^{-1} \left( E_{\al_l}^N T_\al^{-1} v_\la \right)
= [N]_{\al_l}! T_{(w')^{-1}}^{-1} v_\la
\]
and similarly
\[
(\tau E_{\be_l})^m T_{w^{-1}}^{-1} v_\la = 0, 
\; \; \forall \, m >N,
\]
recall \eqref{tau} and \eqref{tau-ident}. Using the formula \eqref{phi1}
for the antihomomorphism 
$\phi_w \colon R^w_0 \to \UU^-[w]$, 
we obtain that for all $\la \in \PP^+$, $\xi \in V(\la)^*$
\begin{align*}
\phi_w( c_\xi^\la e_w^{-\la} ) &= 
\frac{ (q_{\al_l}^{-1} - q_{\al_l})^N}
{q_{\al_l}^{N (N-1)/2}}
\sum_{m_1, \ldots, m_{l-1} \in \Nset}
\left( \prod_{j=1}^{l-1} 
\frac{ (q_{\al_j}^{-1} - q_{\al_j})^{m_j}}
{q_{\al_j}^{m_j (m_j-1)/2} [m_j]_{\al_j}! } \right) 
\\
& \hspace{0.5cm} \times 
\lcor \xi, (\tau E_{\be_1})^{m_1} \ldots 
(\tau E_{\be_{l-1}})^{m_{l-1}} T^{-1}_{(w')^{-1}} v_\la \rcor 
F_{\be_l}^N F_{\be_{l-1}}^{m_{l-1}} \ldots F_{\be_1}^{m_1}
\\
&= 
\frac{ (q_{\al_l}^{-1} - q_{\al_l})^N}
{q_{\al_l}^{N (N-1)/2}}
F_{\be_l}^N \phi_{w'}( c_\xi^\la e_{w'}^{-\la} )
\mod \sum_{m=0}^{N -1 } 
F_{\be_l}^m \UU^-[w'],
\end{align*}
which completes the proof of the proposition.
\qed
\medskip
\\
\noindent
{\em{Proof of \thref{main2-ind} (i).}} In the proof of \thref{main2}
we showed that $l \notin \LP_\ib(y)$ implies 
$F_{\be_l} \notin I_w(y)$. We apply \prref{ind} (i) 
for the iterated Ore extension \eqref{Uwiter} 
and $J = I_w(y)$. Since $I_w(y)$ is a $\Tset^r$-invariant 
completely prime ideal of $\UU^-[w]$, $\lt(I_w(y))$ is a 
$\Tset^r$-invariant completely prime ideal of $\UU^-[w s_{\al_l}]$.
By \thref{Class} (i)
\[
\lt(I_w(y)) = I_{w s_{\al_l}}(y') 
\] 
for some $y' \in W^{\leq w s_\al}$. Let $\la \in \PP^+$ and 
$\xi \in (\UU^- T_y v_\la)^\perp \subset (V(\la))^*$. 
Then $\phi_w( c^\la_\xi e^{-\la}_w) \in I_w(y)$ and 
by \prref{aux2}, 
$\phi_{w s_{\al_l}}( c^\la_\xi e^{-\la}_{w s_{\al_l}}) \in \lt(I_w(y))$. Therefore 
$\lt(I_w(y)) \supseteq I_{w s_{\al_l}}(y)$. Applying 
\thref{Class} (ii), we obtain that $y' \geq y$. 
By \eqref{GK1}, 
\begin{multline*}
\GKdim (\UU^-[w] /I_w(y)) = 
\GKdim (\UU^-[w s_{\al_l} ] /\lt(I_w(y)) ) + 1 
\\
= \GKdim (\UU^-[w s_{\al_l} ] /I_w(y') ) + 1.
\end{multline*}
It follows from \cite[Theorem 5.8]{Y5} that 
\[
\GKdim (\UU^-[w] /I_w(y)) = l - \ell(y) \; \; 
\mbox{and} \; \; 
\GKdim (\UU^-[w s_{\al_l} ] /I_w(y') ) = l- 1 - \ell(y').
\]
Therefore $\ell(y') = \ell(y)$. Since $y' \geq y$, this is only possible 
if $y' = y$, i.e.,
\[
\lt(I_w(y)) = I_{w s_{\al_l}}(y).
\]
\qed
\subsection{Proof of the second part of \thref{main2-ind}}
\label{4.3}
A straightforward computation of the contraction 
$I_w(y) \cap \UU^-[w s_{\al_l}]$ in the Ore extension 
\eqref{Ore-2} is very involved and impractical.
We investigate this contraction in 
a roundabout way by comparing monoids of normal 
elements. We apply \prref{ind} (ii) to deduce 
that 
\begin{equation}
\label{eqq6}
\UU^-[w]/I_w(y)  \cong \UU^-[w s_{\al_l}]/(I_w(y) \cap \UU^-[w s_{\al_l}])
\end{equation}
and \thref{Class} (i) to deduce that 
$I_w(y) \cap \UU^-[w s_{\al_l}] = I_{w s_{\al_l}}(y')$ for some 
$y' \in W^{\leq w s_{\al_l}}$.
From \eqref{commute} we have a supply of nonzero normal elements of the algebras
$\UU^-[w]/I_w(y)$. We prove a characterization of certain 
(equivariantly) normal elements of $\UU^-[w]/I_w(y)$.
With its help we compare the monoids of these equivariantly 
normal elements of the two sides of \eqref{eqq6} and 
deduce that $y' = y s_{\al_l}$.

The weight lattice $\PP$ of $\g$ is embedded in $\Tset^r$ via
$\mu \mt (q^{ \lcor \mu, \al\spcheck \rcor })_{\al \in \Pi}$. The 
$\Tset^r$-action \eqref{torus-act} gives rise to an action of $\PP$ on 
$\UU_q(\g)$, $\UU^-[w]$, and $\UU^-[w]/I_w(y)$, given by 
\[
\mu \cdot x = q^{ \lcor \mu, \ga \rcor} x, 
\; \; \ga \in \QQ, x \in (\UU_q(\g))_\ga.
\]

If a group $M$ acts on a ring $R$ by ring automorphisms, an element 
$u$ of $R$ is called an $M$-normal element if there exists $\mu \in M$
such that 
\[
u x = (\mu \cdot x) u, \; \; \forall \, x \in R. 
\]  
(In relation to equivariant polynormality, in the definition of 
$M$-normal element one sometimes requires that $u$ be an
$M$-eigenvector, see \cite{Y5}. For the sake of clarity, 
we will use the extra term homogeneous to emphasize this.)
Here and below the term homogeneous will refer to the $\QQ$-gradings
of $\UU_q(\g)$, $\UU^-[w]$, and $\UU^-[w]/I_w(y)$. 

By \eqref{commute}, for all $y \in W^{\leq w}$ 
the elements $b_{y,w}^\la$, $\la \in \PP$ are 
nonzero homogeneous $\PP$-normal elements of $\UU^-[w]/I_w(y)$.
The next proposition is a result in the opposite direction 
concerning the possible weights of all
homogeneous $\PP$-normal elements of $\UU^-[w]/I_w(y)$.

\bpr{homog-normal} For all base fields $\KK$, $q \in \KK^*$ not a 
root of unity, Weyl group elements $y \leq w$, and nonzero 
homogeneous $\PP$-normal elements $u \in \UU^-[w]/I_w(y)$,
there exists $\mu \in (1/2)\PP$ such that 
$(w-y) \mu \in \QQ_{\SS(w)}$, $u \in (\UU^-[w]/I_w(y))_{(w-y)\mu}$, 
$(w+y)\mu \in \PP$, and
\[
u x = q^{ - \lcor (w+y)\mu, \ga \rcor } x u, 
\; \; \forall \, \ga \in \QQ, x \in (\UU^-[w]/I_w(y))_\ga.
\]
\epr
\begin{proof} Let $u \in (\UU^-[w]/I_w(y))_{\ga'}$, $\ga' \in \QQ_{\SS(w)}$
be a homogeneous $\PP$-normal element of $\UU^-[w]/I_w(y)$
such that 
\begin{equation}
\label{P-norm}
u x = q^{\lcor \mu' , \ga \rcor } x u, \; \; \forall \, \ga \in \QQ, 
x \in (\UU^-[w]/I_w(y))_\ga
\end{equation}
for some $\mu' \in \PP$.
Eqs. \eqref{commute} and \eqref{P-norm} imply
\[
b^\la_{y,w} u = 
q^{ - \lcor (w+y) \la, \ga' \rcor } u b^\la_{y,w} =
q^{ - \lcor (w+y) \la, \ga' \rcor } 
q^{ \lcor \mu' , (w-y) \la \rcor }
b^\la_{y,w} u
\] 
for all $\la \in \PP^+$. Because $q \in \KK^*$ is not a root of unity
and $\UU^-[w]/I_w(y)$ is a domain
\[
- \lcor (w+y) \la, \ga' \rcor + \lcor \mu', (w-y) \la \rcor = 0, 
\; \; \forall \, \la \in \PP^+.
\] 
Therefore
\[
\lcor w \la, (w y^{-1}+ 1) \ga' \rcor + 
\lcor w \la, (w y^{-1}- 1) \mu' \rcor = 0, \; \; \forall \, 
\la \in \PP^+, 
\]
i.e., 
\[
(w y^{-1} + 1) \ga' = (w y^{-1} -1)( - \mu') = 0.
\]
Using the standard linear algebra argument for Cayley transforms, 
we obtain that there exits $\mu \in \Qset \Pi$ such that 
\begin{equation}
\label{gamu}
\ga' = (w y^{-1} -1) y \mu = (w -y )\mu \; \; 
\mbox{and} \; \; 
- \mu' = (w y^{-1} + 1) y \mu = (w+y )\mu
\end{equation}
(see for instance the proof of \cite[Theorem 3.6]{Y4}). Adding the two equalities 
leads to $2 w(\mu) = \ga' - \mu'$, i.e., 
$\mu = (1/2) w^{-1}(\ga' - \mu') \in (1/2) \PP$. 
Moreover $(w-y)\mu = \ga' \in \QQ_{\SS(w)}$, 
$u \in (\UU^-[w]/I_w(y))_{\ga'} = (\UU^-[w]/I_w(y))_{(w-y)\mu}$, 
and $(w+y)\mu = - \mu' \in \PP$. Finally, substituting \eqref{gamu} 
in \eqref{P-norm} gives
\[
u x = q^{ - \lcor (w+y)\mu, \ga \rcor } x u, \; \; \forall \, \ga \in \QQ, 
x \in (\UU^-[w]/I_w(y))_\ga.
\]
\end{proof}
\noindent
{\em{Proof of \thref{main2-ind} (ii)}}. It was shown in the proof of \thref{main2} 
that $l \in \LP_{\ib}(y)$ implies $F_{\be_l} \in I_w(y)$. 
Recall eq. \eqref{Ore-1}.
Since 
$I_w(y)$ is a $\Tset^r$-invariant completely prime ideal of $\UU^-[w]$, 
$I_w(y) \cap \UU^-[w s_{\al_l}]$ is a $\Tset^r$-invariant 
completely prime ideal of $\UU^-[w s_{\al_l}]$. It follows from 
\thref{Class} (i) that 
\[
I_w(y) \cap \UU^-[w s_{\al_l}] = I_{w s_{\al_l}}(y')
\]
for some $y' \in W^{\leq w s_{\al_l}}$. By \prref{ind} (ii) we have 
the isomorphism of $\QQ$-graded algebras
\[
\UU^-[w s_{\al_l}]/I_{w s_{\al_l}}(y') \cong \UU^-[w]/I_w(y),
\]
because the $\Tset^r$-eigenvectors of $\UU_q(\g)$ with respect to the 
action \eqref{torus-act} are precisely the homogeneous vectors of the 
$\QQ$-grading of $\UU_q(\g)$. Denote the support of 
the $\QQ$-grading of the above algebras:
\[
\QQ' := \Zset \{ \ga \in \QQ \mid (\UU^-[w s_{\al_l}]/I_{w s_{\al_l}}(y'))_\ga \neq 0 \}
\subseteq \QQ.
\]
Let $\la \in \PP$. Eq. \eqref{commute} implies that $b^\la_{y,w}$ is 
a nonzero homogeneous $\PP$-normal element of $\UU^-[w s_{\al_l}]/I_{w s_{\al_l}}(y')$
such that 
\begin{multline}
\label{norm1}
b^\la_{y,w} \in (\UU^-[w s_{\al_l}]/I_{w s_{\al_l}}(y'))_{ (w-y)\la }
\; \; \mbox{and} \\
b_{y,w}^\la x
= 
q^{ - \lcor (w+y)\la, \ga \rcor }
x b_{y,w}^\la, \; \; 
\forall \, 
\ga \in \QQ', 
x \in (\UU^-[w s_{\al_l}]/I_{w s_{\al_l}}(y'))_\ga.
\end{multline}
We apply \prref{homog-normal} to the algebra
$\UU^-[w s_{\al_l}]/I_{w s_{\al_l}}(y')$ and 
the $\PP$-normal element $b^\la_{y,w}$. This shows that 
there exists $\mu' \in (1/2)\PP$ such that 
\begin{multline}
\label{norm2}
b^\la_{y,w} \in (\UU^-[w s_{\al_l}]/I_{w s_{\al_l}}(y'))_{ (w s_{\al_l}-y')\mu },
\\
b_{y,w}^\la x
= 
q^{ - \lcor (w s_{\al_l} +y')\mu, \ga \rcor }
x b_{y,w}^\la, \; \;  
\forall \, 
\ga \in \QQ', 
x \in (\UU^-[w s_{\al_l}]/I_{w s_{\al_l}}(y'))_\ga,
\end{multline}
and $(w s_{\al_l}+y') \mu \in \PP$. Combining \eqref{norm1} and \eqref{norm2}, 
and using that $q \in \KK^*$ is not a root of unity
and $\UU^-[w s_{\al_l}]/I_{w s_{\al_l}}(y')$ is a domain
leads to 
\begin{equation}
\label{twoeqs}
(w-y)\la = (w s_{\al_l}-y')\mu \; \; 
\mbox{and} \; \; 
\lcor (w+y) \la, \ga \rcor =
\lcor (w s_{\al_l} +y')\mu, \ga \rcor, \; \; 
\forall \, \ga \in \QQ'. 
\end{equation}
Therefore 
\begin{multline}
\label{ml}
\lcor w \la, \ga \rcor = 
\lcor (w-y)\la + (w+ y)\la, \ga \rcor 
\\
= 
\lcor (w s_{\al_l} - y')\mu + (w s_{\al_l} +y')\mu, \ga \rcor
= \lcor w s_{\al_l}(\mu), \ga \rcor, \; \; 
\forall \, \ga \in \QQ'.
\end{multline}
For all $\nu \in \PP^+$, $(w s_{\al_l} - y')\nu \in \QQ'$ 
because 
$b^\nu_{w s_{\al_l}, y'} \in 
(\UU^-[w s_{\al_l}]/I_{w s_{\al_l}}(y'))_{(w s_{\al_l} - y')\nu} \backslash \{ 0 \}$.
Hence, by \eqref{ml}
\[
\lcor w s_{\al} ( s_{\al_l} \la - \mu), (w s_{\al_l} - y')\nu \rcor = 0, 
\; \; \forall \, \nu \in \PP^+, 
\]
i.e., 
\[
\lcor (y' - w s_{\al}) ( s_{\al_l} \la - \mu), y' \nu \rcor = 0, 
\; \; \forall \, \nu \in \PP^+.
\]
Thus $(y' - w s_{\al})\mu = (y' - w s_{\al}) s_{\al_l} \la$. 
By taking into account the first part of \eqref{twoeqs}, we obtain
\[
(w-y)\la = ( w s_{\al} - y' ) s_{\al_l} \la.
\] 
Therefore $y \la = y' s_{\al_l} (\la)$ for all $\la \in \PP^+$. We have 
$y' = y s_{\al_l}$ and hence
\[
I_w(y) \cap \UU^-[w s_{\al_l}] = I_{w s_{\al_l}}(y s_{\al_l}),
\]
which completes the proof of part (ii) of \thref{main2}.
\qed

\end{document}